\documentclass[12pt]{article}
\usepackage{amsfonts,amsmath}
\usepackage{epsfig,graphics}

\renewcommand{\appendix}{%
\renewcommand{\section}{%
\newpage
\thispagestyle{plain}%
\secdef\Appendix\sAppendix}%
\setcounter{section}{0}%
\renewcommand{\thesection}{\Alph{section}}%
}
\newcommand{\Appendix}[2][?]{%
\refstepcounter{section}%
\addcontentsline{toc}{Addendum}%
{\protect\numberline{\appendixname~\thesection}#1}%
{\flushleft\LARGE\bfseries\appendixname\ \thesection\par
\centering#2\par}%
\sectionmark{#1}\vspace{\baselineskip}}

\newcommand{\sAppendix}[1]{%
{\flushright\large\bfseries\appendixname\par
\centering#1\par}%
\vspace{\baselineskip}}

\def\be{\begin{equation}}

\def\bea{\begin{eqnarray}}

\def\eea{\end{eqnarray}}

\normalbaselines \oddsidemargin 0.5cm \evensidemargin 0.5cm
\begin{document}

\pagestyle{empty}

\rightline{ CPHT-RR001.01.08}

\vskip 1cm

\begin{center}

{\Huge {\textbf{A New Eight Vertex Model and Higher Dimensional,
Multiparameter Generalizations}}}

\vskip 0.2cm

{\small{\it Dedicated to the memory of Professor Jean Lascoux}}

\vspace{4mm}

{\bf \large B. Abdesselam$^{a,}$\footnote{Email:
boucif@cpht.polytechnique.fr and  boucif@yahoo.fr} and A.
Chakrabarti$^{b,}$\footnote{Email: chakra@cpht.polytechnique.fr}}

\vspace{2mm}

  \emph{$^a$ Laboratoire de Physique Th\'eorique, Universit\'e d'Oran
Es-S\'enia, 31100-Oran, Alg\'erie\\
and \\ Laboratoire de Physique Quantique de la Mati\`ere et de
Mod\'elisations Math\'ematiques, Centre Universitaire de Mascara,
29000-Mascara, Alg\'erie}
  \\
  \vspace{2mm}
  \emph{$^b$ Centre de Physique Th{\'e}orique, Ecole Polytechnique, 91128 Palaiseau Cedex, France.}

\end{center}

\begin{abstract}
{\small \noindent We study statistical models, specifically
transfer matrices corresponding to a multiparameter hierarchy of
braid matrices of $\left(2n\right)^2\times\left(2n\right)^2$
dimensions with $2n^2$ free parameters
$\left(n=1,2,3,\ldots\right)$. The simplest, $4\times 4$ case is
treated in detail. Powerful recursion relations are constructed
giving the dependence on the spectral parameter $\theta$ of the
eigenvalues of the transfer matrix explicitly at each level of
coproduct sequence. A brief study of higher dimensional cases
($n\geq 2$) is presented pointing out features of particular
interest. Spin chain Hamiltonians are also briefly presented for
the hierarchy. In a long final section basic results are
recapitulated with systematic analysis of their contents. Our
eight vertex $4\times 4$ case is compared to standard six vertex
and eight vertex models. }
\end{abstract}

\pagestyle{plain} \setcounter{page}{1}



\section{Introduction}
\setcounter{equation}{0}

In a previous paper \cite{R1} statistical models were presented
starting from a class of multiparameter braid matrices of odd
dimensions (Ref. 1 cites previous sources). This class of
$N^2\times N^2$ braid matrices for $N=2n-1$ ($n=2,3,\ldots$)
depends on $\frac
12\left(N+3\right)\left(N-1\right)=2\left(n^2-1\right)$ parameters
and has a basis of a {\it nested sequence} of projectors
introduced before \cite{R2}. Recently this class has been extended
to include \textit{even} dimensions \cite{R3}. The even
dimensional $\left(2n\right)^2\times\left(2n\right)^2$ braid
matrices depend, apart from the spectral parameter $\theta$, on
$2n^2$ free parameters. One such parameter can, as usual, be
absorbed by a suitable choice of an overall normalization factor.
Such a quadratic increase of free parameters with $n$ (even or
odd) is the most salient, unique feature of our models with their
{\it nested sequence} basis.

We start with a detailed study of the $4\times 4$ case. This is
our {\it new} eight vertex model. It should be compared to {\it
exotic} eight vertex model based on the S\O3 bialgebra \cite{R4}.
There negative Boltzmann weights were involved and were commented
upon. Here, for suitable choices of parameters, one has all
weights positive or zero for $\theta>0$ or for an alternative
choice for $\theta<0$. This $4\times 4$ case has another major
distinguishing feature. It finds its place as the first step in an
explicitly constructed multiparameter hierarchy. This is to be
contrasted with the multistate generalization of six vertex model
(Ref. 5 and the original sources cited there) in which the
parametrization remains restricted to the six vertex level as the
dimensions increases. We will present some remarkable general
features conserved in the entire hierarchy. But here our study
will be brief - a beginning. We hope to explore further elsewhere.
We will also briefly present spin chain Hamiltonians  for our
class of solutions. Some aspects are better discussed after the
constructions are presented. Such comments are reserved for the
end.

\section{Even dimensional multiparameter braid matrices}
\setcounter{equation}{0}

We briefly recapitulate the construction of Ref. 3, closely
related to that for odd dimensions \cite{R1,R2}. The notation of
Ref. 3 are maintained. Define the \textit{nested sequence} of
projectors
\begin{equation}
P_{ij}^{(\epsilon)}=\frac 12\left\{\left(ii\right)\otimes
\left(jj\right)+\left(\bar{i}\bar{i}\right)\otimes
\left(\bar{j}\bar{j}\right)+\epsilon\left[\left(i\bar{i}\right)\otimes
\left(j\bar{j}\right)+\left(\bar{i}i\right)\otimes
\left(\bar{j}j\right)\right]\right\},
\end{equation}
where $i,j\in\left\{1,\cdots,n\right\}$, $\epsilon=\pm$,
$\bar{i}=2n+1-i$, $\bar{j}=2n+1-j$ and $\left(ij\right)$ denotes
the matrix with an unique non-zero element 1 on row $i$ and column
$j$. They provide a complete basis satisfying
\begin{equation}
P_{ij}^{(\epsilon)}P_{kl}^{(\epsilon')}=\delta_{ik}\delta_{jl}\delta_{\epsilon\epsilon'}P_{ij}^{(\epsilon)},\qquad
\hbox{and}\qquad\sum_{\epsilon=\pm}\sum_{i,j=1}^n\left(P_{ij}^{(\epsilon)}+P_{i\bar{j}}^{(\epsilon)}\right)=I_{(2n)^2\times
(2n)^2},
\end{equation}
where $P_{i\bar{j}}^{(\epsilon)}$ is obtained by changing $j$ with
$\bar{j}$ in $P_{ij}^{(\epsilon)}$. The $(2n)^2\times (2n)^2$
braid matrices
\begin{equation}
\hat{R}\left(\theta\right)=\sum_{\epsilon=\pm}\sum_{i,j=1}^ne^{m_{ij}^{(\epsilon)}\theta}\left(P_{ij}^{(\epsilon)}+P_{i\bar{j}}^{(\epsilon)}\right).
\end{equation}
satisfy
\begin{equation}
\widehat{R}_{12}\left(\theta\right)\widehat{R}_{23}\left(\theta+\theta'\right)\widehat{R}_{12}\left(\theta'\right)=
\widehat{R}_{23}\left(\theta'\right)\widehat{R}_{12}\left(\theta+\theta'\right)\widehat{R}_{23}\left(\theta\right),
\end{equation}
where $\widehat{R}_{12}=\widehat{R}\otimes I$ and
$\widehat{R}_{23}=I\otimes \widehat{R}$. Here, as for odd
dimensions \cite{R1,R2} the free parameters enter as exponents
$\left(e^{m_{ij}^{(\epsilon)}\theta}\right)$ and the basic
constraint \cite{R1,R2} $m_{ij}^{(\epsilon)}=
m_{i\bar{j}}^{(\epsilon)}$ has already been incorporated in (2.3).
This leaves $2n^2$ free parameters. The simplest examples are:

\smallskip

\noindent{$\bullet$} \underline{$N=2n=2$}:
\begin{equation}
\hat{R}\left(\theta\right)=\begin{vmatrix}
  a_+ & 0 & 0 & a_- \\
  0 & a_+ & a_- & 0 \\
  0 & a_- & a_+ & 0 \\
  a_- & 0 & 0 & a_+ \\
\end{vmatrix},
\end{equation}
where
\begin{equation}
a_{\pm}=\frac 12\left(e^{m_{11}^{(+)}\theta}\pm
e^{m_{11}^{(-)}\theta}\right).
\end{equation}

\smallskip

\noindent{$\bullet$} \underline{$N=2n=4$}:
\begin{equation}
\hat{R}\left(\theta\right)=\begin{vmatrix}
  D_{11} & 0 & 0 & A_{1\bar{1}} \\
  0 & D_{22} & A_{2\bar{2}} & 0 \\
  0 & A_{\bar{2}2} & D_{\bar{2}\bar{2}} & 0 \\
  A_{\bar{1}1} & 0 & 0 & D_{\bar{1}\bar{1}} \\
\end{vmatrix},
\end{equation}
with
\begin{eqnarray}
&&D_{11}=D_{\bar{1}\bar{1}}=\left(\begin{array}{cccc}
  a_+ & 0 & 0 & 0 \\
  0 & b_+ & 0 & 0 \\
  0 & 0 & b_+ & 0 \\
  0 & 0 & 0 & a_+ \\
\end{array}\right),\qquad
D_{22}=D_{\bar{2}\bar{2}}=\left(\begin{array}{cccc}
  c_+ & 0 & 0 & 0 \\
  0 & d_+ & 0 & 0 \\
  0 & 0 & d_+ & 0 \\
  0 & 0 & 0 & c_+ \\
\end{array}\right),\nonumber\\
&&A_{1\bar{1}}=A_{\bar{1}1}=\left(\begin{array}{cccc}
  0 & 0 & 0 & a_- \\
  0 & 0 & b_- & 0 \\
  0 & b_- & 0 & 0 \\
  a_- & 0 & 0 & 0 \\
\end{array}\right),\qquad
A_{2\bar{2}}=A_{\bar{2}2}=\left(\begin{array}{cccc}
  0 & 0 & 0 & c_- \\
  0 & 0 & d_- & 0 \\
  0 & d_- & 0 & 0 \\
  c_- & 0 & 0 & 0 \\
\end{array}\right)
\end{eqnarray}
and
\begin{eqnarray}
&&a_{\pm}=\frac 12\left(e^{m_{11}^{(+)}\theta}\pm
e^{m_{11}^{(-)}\theta}\right),\qquad b_{\pm}=\frac
12\left(e^{m_{12}^{(+)}\theta}\pm e^{m_{12}^{(-)}\theta}\right),\nonumber\\
&&c_{\pm}=\frac 12\left(e^{m_{21}^{(+)}\theta}\pm
e^{m_{21}^{(-)}\theta}\right),\qquad d_{\pm}=\frac
12\left(e^{m_{22}^{(+)}\theta}\pm
e^{m_{22}^{(-)}\theta}\right)\end{eqnarray} We will often write
\begin{equation}
D_{11}=\left(a_+,b_+,b_+,a_+\right)_{\hbox{diag.}}, \qquad
A_{1\bar{1}}=\left(a_-,b_-,b_-,a_-\right)_{\hbox{anti-diag.}}
\end{equation}
and so on. Generalization for $n>2$ is entirely straightforward.

\section{$\hat{R}\mathbf{TT}$ relations and transfer matrix ($n=1$) }
\setcounter{equation}{0}

To simplify computations, taking out an overall factor $a_+$ we
write $\hat{R}\left(\theta\right)$ of (2.5) as
\begin{equation}
\hat{R}\left(\mathbf{x}\right)=\begin{vmatrix}
  1 & 0 & 0 & \mathbf{x} \\
  0 & 1 & \mathbf{x} & 0 \\
  0 & \mathbf{x} & 1 & 0 \\
  \mathbf{x} & 0 & 0 & 1 \\
\end{vmatrix},
\end{equation}
where
\begin{equation}
\mathbf{x}=\frac{e^{m_{11}^{(+)}\theta}-
e^{m_{11}^{(-)}\theta}}{e^{m_{11}^{(+)}\theta}+
e^{m_{11}^{(-)}\theta}}=\tanh\left(\mu\theta\right),
\end{equation}
with $2\mu=m_{11}^{(+)}-m_{11}^{(-)}$. Note that
\begin{equation}
\left\{m_{11}^{(+)}>m_{11}^{(-)},\,\theta>0\right\}\Longrightarrow
\mathbf{x}>0 \qquad\hbox{and}\qquad
\left\{m_{11}^{(+)}<m_{11}^{(-)},\,\theta<0\right\}\Longrightarrow
\mathbf{x}>0.
\end{equation}
We consider these two domains separately, assuring nonnegative,
real Boltzmann weights. For purely imaginary parameters ${\mathbf
i}m_{11}^{(\pm)}$ ($m_{11}^{(\pm)}$ real)
\begin{equation}
\mathbf{x}={\mathbf i}\tan\left(\mu\theta\right),
\end{equation}
with a normalization factor $\cos\mu\theta$, one obtains a unitary
braid matrix
\begin{equation}
\hat{R}\left(\mathbf{x}\right)^{+}\hat{R}\left(\mathbf{x}\right)=I.
\end{equation}
This has been pointed out in sec. 3 of Ref. 3 citing Ref. 6. We
continue to study here real matrices.

The basic $\hat{R}{\mathbf {TT}}$ equation can now be written as
\begin{equation}
\hat{R}\left(\mathbf{x}''\right)\left({\mathbf
T}\left(\mathbf{x}\right)\otimes {\mathbf
T}\left(\mathbf{x}'\right)\right) =\left({\mathbf
T}\left(\mathbf{x}'\right)\otimes{\mathbf
T}\left(\mathbf{x}\right)\right)\hat{R}\left(\mathbf{x}''\right),
\end{equation}
where
\begin{equation}
\mathbf{x}=\tanh\left(\mu\theta\right),\qquad
\mathbf{x}'=\tanh\left(\mu\theta'\right),\qquad
\mathbf{x}''=\frac{\mathbf{x}-\mathbf{x}'}{1-\mathbf{xx}'}=\tanh\mu\left(\theta-\theta'\right).
\end{equation}
We denote ${\mathbf T}={\mathbf T}\left(\mathbf{x}\right)$,
${\mathbf T}'={\mathbf T}\left(\mathbf{x}'\right)$,
$\hat{R}''=\hat{R} \left(\mathbf{x}''\right)$ and the 4 blocks of
${\mathbf T}$ as
\begin{equation}
{\mathbf T}=\begin{vmatrix}
  {\mathbf T}_{11} & {\mathbf T}_{12} \\
  {\mathbf T}_{21} & {\mathbf T}_{22} \\
  \end{vmatrix}\equiv \begin{vmatrix}
  A & B \\
  C & D \\
  \end{vmatrix}.
\end{equation}
Defining the $2\times 2$ matrices
\begin{equation}
I=\begin{vmatrix}1 & 0 \\
 0 & 1 \\
  \end{vmatrix},\qquad K=\begin{vmatrix}
  0 & 1 \\
  1 & 0 \\
  \end{vmatrix}
\end{equation}
(3.6) is now (with $A=A\left(\mathbf{x}\right)$,
$A'=A\left(\mathbf{x}'\right)$, etc.)
\begin{equation}
\left(I\otimes I+\mathbf{x}''K\otimes K\right)\begin{vmatrix}A & B \\
 C & D \\\end{vmatrix}\otimes\begin{vmatrix}A' & B' \\
 C' & D' \\\end{vmatrix}=\begin{vmatrix}A' & B' \\
 C' & D' \\\end{vmatrix}\otimes\begin{vmatrix}A & B \\
 C & D \\\end{vmatrix}\left(I\otimes I+\mathbf{x}''K\otimes K\right).
\end{equation}
A detailed study of (3.10) is given in Appendix. Recursion
relations will be extracted from it and implemented to construct
the spectrum of the eigenvalues of the transfer matrix in sec. 5.
The transfer matrix $\left(A+D\right)_r\equiv A_r+D_r$ is obtained
by starting with fundamental $2\times 2$ blocks (for $n=1$) and
constructing $2^r\times 2^r$ blocks using standard prescriptions
(coproduct rules - see sec. 5). The starting point, the $2\times
2$ blocks for $r=1$, is provided by the Yang-Baxter matrix
corresponding to (3.1), namely
\begin{equation}
R\left(\mathbf{x}\right)=\mathrm{P}\hat{R}\left(\mathbf{x}\right)=\begin{vmatrix}
  1 & 0 & 0 & \mathbf{x} \\
  0 & \mathbf{x} & 1 & 0 \\
  0 & 1 & \mathbf{x} & 0 \\
  \mathbf{x} & 0 & 0 & 1 \\
\end{vmatrix}\equiv \begin{vmatrix}A_1 & B_1 \\
 C_1 & D_1 \\\end{vmatrix},
\end{equation}
where $\mathrm{P}$ is the standard $4\times 4$ permutation matrix.
Then one uses the coproduct rule
\begin{equation}
\left(\mathbf{T}_{ij}\right)_{r+1}=\sum_k\left(\mathbf{T}_{ik}\right)_{1}\otimes
\left(\mathbf{T}_{kj}\right)_{r}.
\end{equation}
Such a construction guarantees the commutativity of the transfer
matrix , i.e.
\begin{equation}
\left[\mathbf{T}\left(\mathbf{x}\right),\mathbf{T}\left(\mathbf{x}'\right)\right]=0
\end{equation}
quite generally and reducing for our present case to
\begin{equation}
\left[\left(A+D\right)_r,\left(A+D\right)'_r\right]=0
\end{equation}
for all $r$. This is the basic ingredient of exactly solvable
statistical models \cite{R7}. The trace and the highest eigenvalue
of $\mathbf{T}\left(\mathbf{x}\right)$ provide significant
features of the corresponding models. We will obtain these quite
simply and generally for our class. For $n> 1$ certain essential
features will be briefly presented in sec. 6.

\section{Eigenfunctions and eigenvalues of the transfer matrix}
\setcounter{equation}{0}

The next essential step is the construction of eigenstates and
eigenvalues of the transfer matrix. In this section we display the
results for $r=1,2,3,4$. They provide explicit examples of the
iterative structure of the transfer matrices $\mathbf{T}_r$ to be
derived in the next section. Moreover they illustrate how
multiplets involving $r$-th roots of unity and possible
multiplicities of them combine to provide a complete basis of
mutually orthogonal eigenstates of $\mathbf{T}_r$ spanning the
base space. Referring back to them one grasps better the full
content of the general formalism of sec. 5. We start with the
notations of sec. 3.

\smallskip

\noindent{$\bullet$ \underline{$r=1$}:}
\begin{equation}
\mathbf{T}_1=\left(A+D\right)_1=\left(1+\mathbf{x}\right)\begin{vmatrix}
  1 & 0 \\
  0 & 1 \\
  \end{vmatrix},\qquad \left(A-D\right)_1=\left(1-\mathbf{x}\right)\begin{vmatrix}
  1 & 0 \\
  0 & -1 \\
  \end{vmatrix}.\end{equation}
The eigenstates are evidently
\begin{equation}
\left|1\right\rangle\equiv\left|\begin{matrix}
  1\\
  0\\
  \end{matrix}\right\rangle,\qquad
\left|\overline{1}\right\rangle\equiv\left|\begin{matrix}
  0\\
  1\\
  \end{matrix}\right\rangle\end{equation}
with
\begin{equation}
\mathbf{T}_1\left(\left|1\right\rangle,
\left|\overline{1}\right\rangle\right) =\left(1+\mathbf{x}\right)
\left(\left|1\right\rangle, \left|\overline{1}\right\rangle\right)
\end{equation} and
\begin{equation}
\hbox{Tr}\left(\mathbf{T}_1\right)=2\left(1+\mathbf{x}\right).
\end{equation}

\smallskip

\noindent{$\bullet$ \underline{$r=2$}:}
\begin{equation}
\mathbf{T}_2=\left(A+D\right)_2=\left(1+\mathbf{x}\right)^2\mathbf{X}_{(2,0)}+(1-\mathbf{x})^2\mathbf{X}_{(0,2)},
\end{equation}
where
\begin{equation}
\mathbf{X}_{(2,0)}=\frac 12\begin{vmatrix}
  I_{(1)} & K_{(1)}\\
  K_{(1)} & I_{(1)}\\
  \end{vmatrix}\otimes \begin{vmatrix}
  1 & 0\\
  0 & 1\\
  \end{vmatrix},\qquad \mathbf{X}_{(0,2)}=\frac 12\begin{vmatrix}
  I_{(1)} & K_{(1)}\\
  -K_{(1)} & -I_{(1)}\\
  \end{vmatrix}\otimes \begin{vmatrix}
  1 & 0\\
  0 & -1\\
  \end{vmatrix},
\end{equation}
with $I_{(1)}=\begin{vmatrix}
  1 & 0\\
  0 & 1\\
  \end{vmatrix}$ and $K_{(1)}=\begin{vmatrix}
  0 & 1\\
  1 & 0\\
  \end{vmatrix}$. One has
\begin{equation}
\mathbf{X}_{(2,0)}\mathbf{X}_{(0,2)}=
\mathbf{X}_{(0,2)}\mathbf{X}_{(2,0)}=0.
\end{equation}
We have anticipated the iterative structure of sec. 5. That
$\mathbf{T}_2$ obtained by a straightforward use of coproduct
rules as
\begin{equation}
\mathbf{T}_2=\begin{vmatrix}
1+\mathbf{x}^2 & 0 & 0 & 2\mathbf{x}\\
0 & 2\mathbf{x} & \left(1+\mathbf{x}\right)^2 & 0\\
0 & \left(1+\mathbf{x}\right)^2 & 2\mathbf{x} & 0\\
2\mathbf{x} & 0 & 0 & 1+\mathbf{x}^2\\
\end{vmatrix}
\end{equation}
can be expressed as (4.5) on a basis satisfying (4.7) is the
central lesson. Denote
\begin{equation}
\left|\begin{matrix}
  1\\
  0\\
  \end{matrix}\right\rangle\otimes \left|\begin{matrix}
  1\\
  0\\
  \end{matrix}\right\rangle\equiv\left|11\right\rangle,\qquad
\left|\begin{matrix}
  0\\
  1\\
  \end{matrix}\right\rangle\otimes \left|\begin{matrix}
  0\\
  1\\
  \end{matrix}\right\rangle\equiv\left|\bar{1}\bar{1}\right\rangle,\qquad
\left|\begin{matrix}
  1\\
  0\\
  \end{matrix}\right\rangle\otimes \left|\begin{matrix}
  0\\
  1\\
  \end{matrix}\right\rangle\equiv\left|1\bar{1}\right\rangle,\qquad
\left|\begin{matrix}
  0\\
  1\\
  \end{matrix}\right\rangle\otimes \left|\begin{matrix}
  1\\
  0\\
  \end{matrix}\right\rangle\equiv\left|\bar{1}1\right\rangle.
  \end{equation}
For $r>1$ the order of the indices $\left(1,\bar{1}\right)$
indicates the structure of the tensor product. For example,
$\left|1\bar{1}1\right\rangle= \left|\begin{matrix}
  1\\
  0\\
  \end{matrix}\right\rangle\otimes \left|\begin{matrix}
  0\\
  1\\
  \end{matrix}\right\rangle\otimes
\left|\begin{matrix}
  1\\
  0\\
  \end{matrix}\right\rangle$. One obtains, with upper or lower signs,
\begin{equation}
\mathbf{T}_2\left(\left|11\right\rangle\pm\left|\bar{1}\bar{1}\right\rangle\right)
=\left(1\pm\mathbf{x}\right)^2\left(\left|11\right\rangle\pm\left|\bar{1}\bar{1}\right\rangle\right),
\qquad\mathbf{T}_2\left(\left|1\bar{1}\right\rangle\pm\left|\bar{1}1\right\rangle\right)
=\pm\left(1\pm\mathbf{x}\right)^2\left(\left|1\bar{1}\right\rangle\pm\left|\bar{1}1\right\rangle\right)
\end{equation}
and
\begin{equation}
\hbox{Tr}\left(\mathbf{T}_2\right)=2\left(1+ \mathbf{x}\right)^2.
\end{equation}
The eigenfunctions of $\mathbf{X}_{(2,0)}$ (resp.
$\mathbf{X}_{(0,2)}$) are annihilated by $\mathbf{X}_{(0,2)}$
(resp. $\mathbf{X}_{(2,0)}$) consistently with (4.7). Finally,
again anticipating sec. 5,
\begin{equation}
\left(A-D\right)_2=\frac 12\left(1+\mathbf{x}\right)
\begin{vmatrix}
  I_{(1)} & -K_{(1)}\\
  -K_{(1)} & I_{(1)}\\
  \end{vmatrix}\otimes
\left(A-D\right)_1 +
 \frac 12\left(1-\mathbf{x}\right)\begin{vmatrix}
  I_{(1)} & -K_{(1)}\\
  K_{(1)} & -I_{(1)}\\
  \end{vmatrix}\otimes \left(A+D\right)_1.
\end{equation}

\smallskip

\noindent{$\bullet$ \underline{$r=3$}:} From this stage onwards
one can better appreciate the role of multiplets involving roots
of unity in constructing eigenstates. One obtains by implementing
coproducts in step $(r=2)\longrightarrow(r=3)$,
\begin{equation}
\mathbf{T}_3=\left(A+D\right)_3=\left(1+\mathbf{x}\right)^3\mathbf{X}_{(3,0)}+(1+\mathbf{x})(1-\mathbf{x})^2\mathbf{X}_{(1,2)},
\end{equation}
where the $8\times 8$ matrices satisfy
\begin{equation}
\mathbf{X}_{(3,0)}\mathbf{X}_{(1,2)}=\mathbf{X}_{(1,2)}\mathbf{X}_{(3,0)}=0
\end{equation}
and they are $\mathbf{x}$-independent. They are easy to obtain but
will not, for brevity, be presented explicitly. The base space
splits up into 4-dim. ones, closed under the action of
$\mathbf{T}_3$ and are characterized by even (odd) number of the
index 1 (considering zero as even) respectively. Define
\begin{eqnarray}
&&\left|e_1\right\rangle=
\left|\bar{1}\bar{1}\bar{1}\right\rangle+\left|\bar{1}11\right\rangle+\left|1\bar{1}1\right\rangle+\left|11\bar{1}\right\rangle,\nonumber\\
&&\left|e_2\right\rangle=
-3\left|\bar{1}\bar{1}\bar{1}\right\rangle+\left|\bar{1}11\right\rangle+\left|1\bar{1}1\right\rangle+\left|11\bar{1}\right\rangle,\nonumber\\
&&\left|e_3\right\rangle=\left|\bar{1}11\right\rangle+\omega \left|1\bar{1}1\right\rangle+\omega^2\left|11\bar{1}\right\rangle,\nonumber\\
&&\left|e_4\right\rangle=\left|\bar{1}11\right\rangle+\omega^2\left|1\bar{1}1\right\rangle+\omega\left|11\bar{1}\right\rangle,\end{eqnarray}
where $\omega=e^{\mathbf{i}\frac{2\pi}3}$
$\left(\omega^3=1\,\hbox{and}\,1+\omega+\omega^2=0\right)$. One
has $\left\langle e_i|e_j\right\rangle=0$, $i\neq j$, where
$\left\langle e_i\right|$ denotes the transform of
$\left|e_i\right\rangle$ with conjugated coefficients. Denote also
\begin{eqnarray}
&&\left|o_1\right\rangle=
\left|111\right\rangle+\left|1\bar{1}\bar{1}\right\rangle+\left|\bar{1}1\bar{1}\right\rangle+\left|\bar{1}\bar{1}1\right\rangle,
\nonumber\\
&&\left|o_2\right\rangle=
-3\left|111\right\rangle+\left|1\bar{1}\bar{1}\right\rangle+\left|\bar{1}1\bar{1}\right\rangle+\left|\bar{1}\bar{1}1\right\rangle,\nonumber\\
&&\left|o_3\right\rangle=\left|1\bar{1}\bar{1}\right\rangle+\omega \left|\bar{1}1\bar{1}\right\rangle+\omega^2\left|\bar{1}\bar{1}1\right\rangle,
\nonumber\\
&&\left|o_4\right\rangle=\left|1\bar{1}\bar{1}\right\rangle+\omega^2\left|\bar{1}1\bar{1}\right\rangle+\omega\left|\bar{1}\bar{1}1\right\rangle.
\end{eqnarray}
The states
$\left\{\left|e_i\right\rangle,\left|o_j\right\rangle\right\}$
form a complete basis of orthogonal states. Consistently with
(4.13), (4.14) one obtains
\begin{eqnarray}
&&\mathbf{T}_3\left|e_1\right\rangle=\left(1+\mathbf{x}\right)^3\left|e_1\right\rangle,
\nonumber\\
&&\mathbf{T}_3\left(\left|e_2\right\rangle,\left|e_3\right\rangle,\left|e_4\right\rangle\right)=\left(1+\mathbf{x}\right)\left(1-\mathbf{x}\right)^2
\left(\left|e_2\right\rangle,\omega\left|e_3\right\rangle,\omega^2\left|e_4\right\rangle\right),\nonumber\\
&&\mathbf{T}_3\left|o_1\right\rangle=\left(1+\mathbf{x}\right)^3\left|o_1\right\rangle,\nonumber\\
&&\mathbf{T}_3\left(\left|o_2\right\rangle,\left|o_3\right\rangle,\left|o_4\right\rangle\right)=\left(1+\mathbf{x}\right)\left(1-\mathbf{x}\right)^2
\left(\left|o_2\right\rangle,\omega\left|o_3\right\rangle,\omega^2\left|o_4\right\rangle\right)
\end{eqnarray}
and, finally,
\begin{eqnarray}
&&\hbox{Tr}\left(\mathbf{T}_3\right)=2\left(1+\mathbf{x}\right)^3+2\left(1+\mathbf{x}\right)\left(1-\mathbf{x}\right)^2\left(1+\omega+\omega^2\right)
=2\left(1+\mathbf{x}\right)^3.
\end{eqnarray}

\smallskip

\noindent{$\bullet$ \underline{$r=4$}:} Now
\begin{equation}
\mathbf{T}_4=\left(A+D\right)_4=\left(1+\mathbf{x}\right)^4\mathbf{X}_{(4,0)}+\left(1+\mathbf{x}\right)^2\left(1-\mathbf{x}\right)^2\mathbf{X}_{(2,2)}+
\left(1-\mathbf{x}\right)^4\mathbf{X}_{(0,4)},
\end{equation}
where the $16\times 16$ constant matrices satisfy
\begin{equation}
\mathbf{X}_{(4,0)}\mathbf{X}_{(2,2)}=\mathbf{X}_{(4,0)}\mathbf{X}_{(0,4)}=
\mathbf{X}_{(2,2)}\mathbf{X}_{(0,4)}=\mathbf{X}_{(2,2)}\mathbf{X}_{(4,0)}=
\mathbf{X}_{(0,4)}\mathbf{X}_{(4,0)}=\mathbf{X}_{(0,4)}\mathbf{X}_{(2,2)}=0
\end{equation}
as consequence of recursion relations involved in $\left(A\pm
D\right)_3\longrightarrow \left(A\pm D\right)_4$. The
$\mathbf{X}$'s are obtained fairly easily. Now the even and odd
subspaces are 8-dim. We display the eigenstates explicitly to
illustrate a new feature. Now $r=4$ is not a prime number and 4-th
roots and square roots of unity (corresponding to $r=2\times 2$)
both contribute multiplets involving respective coefficients
$\left(1,\mathbf{i},-1,-\mathbf{i}\right)$ and $\left(1,-1\right)$
(For $r=6$ one would have thus 2-plets, 3-plets and 6-plets). The
$\left|e\right\rangle$ and $\left|o\right\rangle$ spaces have the
orthogonal bases
\begin{eqnarray}
&&\left|e_1\right\rangle=
\left|1111\right\rangle+\left|11\bar{1}\bar{1}\right\rangle+\left|1\bar{1}1\bar{1}\right\rangle+\left|1\bar{1}\bar{1}1\right\rangle
+\left|\bar{1}11\bar{1}\right\rangle+\left|\bar{1}1\bar{1}1\right\rangle+\left|\bar{1}\bar{1}11\right\rangle+
\left|\bar{1}\bar{1}\bar{1}\bar{1}\right\rangle,\nonumber\\
&&\left|e_2\right\rangle=
\left|1111\right\rangle-\left|11\bar{1}\bar{1}\right\rangle+\left|1\bar{1}1\bar{1}\right\rangle-\left|1\bar{1}\bar{1}1\right\rangle
-\left|\bar{1}11\bar{1}\right\rangle+\left|\bar{1}1\bar{1}1\right\rangle-\left|\bar{1}\bar{1}11\right\rangle+
\left|\bar{1}\bar{1}\bar{1}\bar{1}\right\rangle,\nonumber\\
&&\left|e_3\right\rangle= \left|1111\right\rangle-
\left|\bar{1}\bar{1}\bar{1}\bar{1}\right\rangle,\nonumber\\
&&\left|e_4\right\rangle=
\left|1111\right\rangle-\left|1\bar{1}1\bar{1}\right\rangle-\left|\bar{1}1\bar{1}1\right\rangle+
\left|\bar{1}\bar{1}\bar{1}\bar{1}\right\rangle,\nonumber\\
&&\left|e_5\right\rangle=
\left|11\bar{1}\bar{1}\right\rangle-\left|1\bar{1}\bar{1}1\right\rangle-\left|\bar{1}11\bar{1}\right\rangle+
\left|\bar{1}\bar{1}11\right\rangle,\nonumber\\
&&\left|e_6\right\rangle=
\left|1\bar{1}1\bar{1}\right\rangle-\left|\bar{1}1\bar{1}1\right\rangle,\nonumber\\
&&\left|e_7\right\rangle=
\left|11\bar{1}\bar{1}\right\rangle-\mathbf{i}\left|1\bar{1}\bar{1}1\right\rangle+\mathbf{i}\left|\bar{1}11\bar{1}\right\rangle-
\left|\bar{1}\bar{1}11\right\rangle,\nonumber\\
&&\left|e_8\right\rangle=
\left|11\bar{1}\bar{1}\right\rangle+\mathbf{i}\left|1\bar{1}\bar{1}1\right\rangle-\mathbf{i}\left|\bar{1}11\bar{1}\right\rangle-
\left|\bar{1}\bar{1}11\right\rangle,
\end{eqnarray}
and
\begin{eqnarray}
&&\left|o_1\right\rangle=
\left|111\bar{1}\right\rangle+\left|11\bar{1}1\right\rangle+\left|1\bar{1}11\right\rangle+\left|\bar{1}111\right\rangle
+\left|1\bar{1}\bar{1}\bar{1}\right\rangle+\left|\bar{1}1\bar{1}\bar{1}\right\rangle+\left|\bar{1}\bar{1}1\bar{1}\right\rangle+
\left|\bar{1}\bar{1}\bar{1}1\right\rangle,\nonumber\\
&&\left|o_2\right\rangle=
\left|111\bar{1}\right\rangle-\left|11\bar{1}1\right\rangle+\left|1\bar{1}11\right\rangle-\left|\bar{1}111\right\rangle
+\left|\bar{1}\bar{1}\bar{1}1\right\rangle-\left|\bar{1}\bar{1}1\bar{1}\right\rangle+\left|\bar{1}1\bar{1}\bar{1}\right\rangle-
\left|1\bar{1}\bar{1}\bar{1}\right\rangle,\nonumber\\
&&\left|o_3\right\rangle=
\left|111\bar{1}\right\rangle+\left|11\bar{1}1\right\rangle+\left|1\bar{1}11\right\rangle+\left|\bar{1}111\right\rangle
-\left|\bar{1}\bar{1}\bar{1}1\right\rangle-\left|\bar{1}\bar{1}1\bar{1}\right\rangle-\left|\bar{1}1\bar{1}\bar{1}\right\rangle-
\left|1\bar{1}\bar{1}\bar{1}\right\rangle,\nonumber\\
&&\left|o_4\right\rangle=
\left|111\bar{1}\right\rangle-\left|11\bar{1}1\right\rangle+\left|1\bar{1}11\right\rangle-\left|\bar{1}111\right\rangle
-\left|\bar{1}\bar{1}\bar{1}1\right\rangle+\left|\bar{1}\bar{1}1\bar{1}\right\rangle-\left|\bar{1}1\bar{1}\bar{1}\right\rangle+
\left|1\bar{1}\bar{1}\bar{1}\right\rangle,\nonumber\\
&&\left|o_5\right\rangle=
\left|111\bar{1}\right\rangle-\mathbf{i}\left|11\bar{1}1\right\rangle-\left|1\bar{1}11\right\rangle+\mathbf{i}
\left|\bar{1}111\right\rangle,\nonumber\\
&&\left|o_6\right\rangle=
\left|1\bar{1}\bar{1}\bar{1}\right\rangle+\mathbf{i}\left|\bar{1}1\bar{1}\bar{1}\right\rangle-\left|\bar{1}\bar{1}1\bar{1}\right\rangle-\mathbf{i}
\left|\bar{1}\bar{1}\bar{1}1\right\rangle,\nonumber\\
&&\left|o_7\right\rangle=
\left|111\bar{1}\right\rangle+\mathbf{i}\left|11\bar{1}1\right\rangle-\left|1\bar{1}11\right\rangle-
\mathbf{i}\left|\bar{1}111\right\rangle,\nonumber\\
&&\left|o_8\right\rangle=
\left|1\bar{1}\bar{1}\bar{1}\right\rangle-\mathbf{i}\left|\bar{1}1\bar{1}\bar{1}\right\rangle-\left|\bar{1}\bar{1}1\bar{1}\right\rangle+\mathbf{i}
\left|\bar{1}\bar{1}\bar{1}1\right\rangle.
\end{eqnarray}
Define
\begin{equation}
\mathbf{T}_4\left|e_i\right\rangle=\upsilon_i^{(e)}\left|e_i\right\rangle,\qquad
\mathbf{T}_4\left|o_i\right\rangle=\upsilon_i^{(o)}\left|o_i\right\rangle,\qquad
i=1,\ldots, 8.
\end{equation}
Then one obtains, in order,
\begin{eqnarray}
&&\upsilon_1^{(e)},\ldots,\upsilon_8^{(e)}=\left(1+\mathbf{x}\right)^4,\left(1-\mathbf{x}\right)^4,
\left(1-\mathbf{x}\right)^2\left(1+\mathbf{x}\right)^2\left(1,1,-1,-1,\mathbf{i},-\mathbf{i}\right),\nonumber\\
&&\upsilon_1^{(o)},\ldots,\upsilon_8^{(o)}=\left(1+\mathbf{x}\right)^4,-\left(1-\mathbf{x}\right)^4,
\left(1-\mathbf{x}\right)^2\left(1+\mathbf{x}\right)^2\left(1,-1,\mathbf{i},\mathbf{i},-\mathbf{i},-\mathbf{i}\right)
\end{eqnarray}
and
\begin{equation}
\hbox{Tr}\left(\mathbf{T}_4\right)=2\left(1+\mathbf{x}\right)^4.
\end{equation}

\section{Relating $\left(A,B,C,D\right)$ for all $r$ and constructing iteratively the eigenvalue spectrum}
\setcounter{equation}{0}

We consider below exclusively the $4\times 4$ case. Coproduct
rules lead, for $N=2$, to the recursion relations
\begin{eqnarray}
&&A_{r+1}=A_1\otimes A_r+B_1\otimes C_r,\qquad D_{r+1}=D_1\otimes
D_r+C_1\otimes B_r,\nonumber\\
&&B_{r+1}=A_1\otimes B_r+B_1\otimes D_r,\qquad C_{r+1}=C_1\otimes
A_r+D_1\otimes C_r,
\end{eqnarray}
where
\begin{equation}
A_1=\begin{vmatrix}
  1 & 0 \\
  0 & \mathbf{x} \\
\end{vmatrix}, \qquad D_1=\begin{vmatrix}
  \mathbf{x} & 0 \\
  0 & 1 \\
\end{vmatrix},\qquad B_1=\begin{vmatrix}
  0 & \mathbf{x} \\
  1 & 0 \\
\end{vmatrix},\qquad C_1=\begin{vmatrix}
  0 & 1 \\
  \mathbf{x} & 0 \\
\end{vmatrix}.
\end{equation}
Thus
\begin{equation}
A_{r+1}=\begin{vmatrix}
  A_r & \mathbf{x}C_r \\
  C_r & \mathbf{x}A_r \\
\end{vmatrix}, \,\,\,
D_{r+1}=\begin{vmatrix}
  \mathbf{x}D_r & B_r \\
  \mathbf{x}B_r & D_r \\
\end{vmatrix},\,\,\,
B_{r+1}=\begin{vmatrix}
  B_r & \mathbf{x}D_r \\
  D_r & \mathbf{x}B_r \\
\end{vmatrix},\,\,\,
C_{r+1}=\begin{vmatrix}
  \mathbf{x}C_r & A_r \\
  \mathbf{x}A_r & C_r \\
\end{vmatrix}.
\end{equation}
Denote $I\equiv I_2=\begin{vmatrix}
  1 & 0 \\
  0 & 1 \\
\end{vmatrix}$, $K\equiv K_2=\begin{vmatrix}
  0 & 1 \\
  1 & 0 \\
\end{vmatrix}$ and for $p$ factors $I_{(p)}=I\otimes I\otimes
\ldots\otimes I$, $K^{(p)}=K\otimes K\otimes \ldots\otimes K$.
Starting from $B_1=KA_1$, $C_1=KD_1$, $D_1=KA_1K$ and iterating
one obtains
\begin{eqnarray}
&&B_r=\left(I_{(r-1)}\otimes K\right)A_r,\qquad
C_r=\left(I_{(r-1)}\otimes K\right)D_r,\nonumber\\
&&D_r=K^{(r)}A_rK^{(r)},\qquad A_r=K^{(r)}D_rK^{(r)}.
\end{eqnarray}
Using $K^{(r)}\left(I_{(r-1)}\otimes K\right)=K^{(r-1)}\otimes I$
one can express $\left(A,B,C,D\right)_{r}$ in term of any one of
them. In particular,
\begin{equation}
\left(B_r\pm C_{r}\right)=\left(I_{(r-1)}\otimes
K\right)\left(A_r\pm D_{r}\right),
\end{equation}
where
\begin{equation}
\left(I_{(r-1)}\otimes
K\right)=\left(K,K,\ldots,K\right)_{\hbox{diag.}} \equiv K_{(r)}
\end{equation}
Thus
\begin{equation}
\left(B_3\pm C_3\right)=\begin{vmatrix}
  K & 0 & 0 & 0 \\
  0 & K & 0 & 0 \\
  0 & 0 & K & 0 \\
  0 & 0 & 0 & K \\
\end{vmatrix}\left(A_3\pm D_3\right).
\end{equation}
Thus exchanging members of successive pairs of rows
$\left[(1,2),(3,4),\ldots,(2p-1,2p),\ldots\right]_{\hbox{(rows)}}$
one obtains $\left(B_r\pm C_r\right)$ from $\left(A_r\pm
D_r\right)$. The following recursion relations are implied (with
the definition (5.6) of $K_{(r)}$)
\begin{eqnarray}
&&\left(A+D\right)_{r+1}=\begin{vmatrix}
  \left(A+\mathbf{x}D\right)_r & \left(B+\mathbf{x}C\right)_r \\
  \left(\mathbf{x}B+C\right)_r & \left(\mathbf{x}A+D\right)_r \\
\end{vmatrix}
\nonumber\\&&\phantom{\left(A+D\right)_{r+1}}=\frac
12\left(1+\mathbf{x}\right)\begin{vmatrix}
  I_{(r)} & K_{(r)} \\
  K_{(r)} & I_{(r)} \\
\end{vmatrix}\begin{vmatrix}
  \left(A+D\right)_r & 0 \\
  0 & \left(A+D\right)_r \\
\end{vmatrix}\nonumber\\
&&\phantom{\left(A+D\right)_{r+1}} + \frac
12\left(1-\mathbf{x}\right)\begin{vmatrix}
  I_{(r)} & K_{(r)} \\
  -K_{(r)} & -I_{(r)} \\
\end{vmatrix}\begin{vmatrix}
  \left(A-D\right)_r & 0 \\
  0 & \left(A-D\right)_r \\
\end{vmatrix}\\
&&\left(A-D\right)_{r+1}=\frac
12\left(1+\mathbf{x}\right)\begin{vmatrix}
  I_{(r)} & -K_{(r)} \\
  -K_{(r)} & I_{(r)} \\
\end{vmatrix}\begin{vmatrix}
  \left(A-D\right)_r & 0 \\
  0 & \left(A-D\right)_r \\
\end{vmatrix}\nonumber\\
&&\phantom{\left(A-D\right)_{r+1}} + \frac
12\left(1-\mathbf{x}\right)\begin{vmatrix}
  I_{(r)} & -K_{(r)} \\
  K_{(r)} & -I_{(r)} \\
\end{vmatrix}\begin{vmatrix}
  \left(A+D\right)_r & 0 \\
  0 & \left(A+D\right)_r \\
\end{vmatrix}
\end{eqnarray}
The signification of these relations concerning eigenvalues and
why we display also $\left(A-D\right)_{r+1}$ will be explained
below.

Let us introduce at this point the general possibilities
\begin{eqnarray}
&&A_{r_1+r_2}=A_{r_1}\otimes A_{r_2}+B_{r_1}\otimes C_{r_2},\qquad
D_{r_1+r_2}=D_{r_1}\otimes D_{r_2}+C_{r_1}\otimes
B_{r_2},\nonumber\\
&&B_{r_1+r_2}=A_{r_1}\otimes B_{r_2}+B_{r_1}\otimes D_{r_2},\qquad
C_{r_1+r_2}=C_{r_1}\otimes A_{r_2}+D_{r_1}\otimes C_{r_2}.
\end{eqnarray}
Since the sequence for odd and even $r$ have some distinct typical
features a {\it two-step} iteration can be of interest for
$\left(A\pm D\right)$. One has (in evident notations)
\begin{eqnarray}
&&2\left(A+D\right)_{r+2}=\left(A+D\right)_{2}\otimes
\left(A+D\right)_{r}+\left(A-D\right)_{2}\otimes
\left(A-D\right)_{r}\nonumber\\
&&\phantom{2\left(A+D\right)_{r+2}}+\left(B+C\right)_{2}\otimes
\left(B+C\right)_{r}-\left(B-C\right)_{2}\otimes
\left(B-C\right)_{r},\\
&&2\left(A-D\right)_{r+2}=\left(A+D\right)_{2}\otimes
\left(A-D\right)_{r}+\left(A-D\right)_{2}\otimes
\left(A+D\right)_{r}\nonumber\\
&&\phantom{2\left(A+D\right)_{r+2}}-\left(B+C\right)_{2}\otimes
\left(B-C\right)_{r}+\left(B-C\right)_{2}\otimes
\left(B+C\right)_{r}
\end{eqnarray}
leading to
\begin{eqnarray}
&&4\left(A+D\right)_{r+2}=\nonumber\\
&&\left(1+\mathbf{x}\right)^2
\begin{vmatrix}
  I_{(r)} & K_{(r)} & K_{(r)} & I_{(r)} \\
  K_{(r)} & I_{(r)} & I_{(r)} & K_{(r)} \\
  K_{(r)} & I_{(r)} & I_{(r)} & K_{(r)} \\
  I_{(r)} & K_{(r)} & K_{(r)} & I_{(r)} \\
\end{vmatrix}
\begin{vmatrix}
  \left(A+D\right)_r & 0 & 0 & 0 \\
  0 & \left(A+D\right)_r & 0 & 0 \\
  0 & 0 & \left(A+D\right)_r & 0 \\
  0 & 0 & 0 & \left(A+D\right)_r \\
\end{vmatrix}\\
&&+\left(1-\mathbf{x}\right)^2
\begin{vmatrix}
  I_{(r)} & -K_{(r)} & K_{(r)} & -I_{(r)} \\
  K_{(r)} & -I_{(r)} & I_{(r)} & -K_{(r)} \\
  -K_{(r)} & I_{(r)} & -I_{(r)} & K_{(r)} \\
  -I_{(r)} & K_{(r)} & -K_{(r)} & I_{(r)} \\
\end{vmatrix}
\begin{vmatrix}
  \left(A+D\right)_r & 0 & 0 & 0 \\
  0 & \left(A+D\right)_r & 0 & 0 \\
  0 & 0 & \left(A+D\right)_r & 0 \\
  0 & 0 & 0 & \left(A+D\right)_r \\
\end{vmatrix}\nonumber\\
&&+\left(1+\mathbf{x}\right)\left(1-\mathbf{x}\right)
\begin{vmatrix}
  I_{(r)} & 0 & K_{(r)} & 0 \\
  -K_{(r)} & 0 & -I_{(r)} & 0 \\
  0 & I_{(r)} & 0 & K_{(r)} \\
  0 & -K_{(r)} & 0 & -I_{(r)} \\
\end{vmatrix}
\begin{vmatrix}
  \left(A-D\right)_r & 0 & 0 & 0 \\
  0 & \left(A-D\right)_r & 0 & 0 \\
  0 & 0 & \left(A-D\right)_r & 0 \\
  0 & 0 & 0 & \left(A-D\right)_r \\
\end{vmatrix}\nonumber
\end{eqnarray}
From (5.12) one obtains an analogous result for
$\left(A-D\right)_{r+2}$.

We now extract one fundamental consequence of the recursion
relations (5.8-13). Starting with
$\left(A+D\right)_{1}=\left(1+\mathbf{x}\right)\begin{vmatrix}
  1 & 0 \\
  0 & 1 \\
\end{vmatrix}$ and $\left(A-D\right)_{1}=\left(1-\mathbf{x}\right)\begin{vmatrix}
  1 & 0 \\
  0 & -1 \\
\end{vmatrix}$ these relations imply (with $\mathbf{x}$-independent
matrices $\mathbf{X}$)
\begin{eqnarray}
&&\left(A+D\right)_{r}=\left(1+\mathbf{x}\right)^r\mathbf{X}_{(r,0)}+\left(1+\mathbf{x}\right)^{r-2}\left(1-\mathbf{x}\right)^{2}\mathbf{X}_{(r-2,2)}
\nonumber\\
&&\phantom{\left(A+D\right)_{r}=}+\cdots+\left(1+\mathbf{x}\right)^2\left(1-\mathbf{x}\right)^{r-2}\mathbf{X}_{(2,r-2)}+
\left(1-\mathbf{x}\right)^{r}\mathbf{X}_{(0,r)}
\end{eqnarray}
for even $r$ and
\begin{eqnarray}
&&\left(A+D\right)_{r}=\left(1+\mathbf{x}\right)^r\mathbf{X}_{(r,0)}+\cdots+\left(1+\mathbf{x}\right)\left(1-\mathbf{x}\right)^{r-1}\mathbf{X}_{(1,r-1)}
\end{eqnarray}
for odd $r$, powers of $\left(1-\mathbf{x}\right)$ being always
even in both cases. Correspondingly (again with
$\mathbf{x}$-independent $2^r\times 2^r$ matrices $\mathbf{Y}$)
\begin{eqnarray}
&&\left(A-D\right)_{r}=\left(1+\mathbf{x}\right)^{r-1}\left(1-\mathbf{x}\right)\mathbf{Y}_{(r-1,1)}+
\left(1+\mathbf{x}\right)^{r-3}\left(1-\mathbf{x}\right)^3\mathbf{Y}_{(r-3,3)}+\cdots+\nonumber\\
&&\phantom{\left(A-D\right)_{r}=}\left(1+\mathbf{x}\right)^\delta
\left(1-\mathbf{x}\right)^{r-\delta}\mathbf{Y}_{(\delta,r-\delta)},
\end{eqnarray}
where $\delta=\frac 12\left(1+\left(-1\right)^r\right)$. There
being here only odd powers of $\left(1-\mathbf{x}\right)$. When
these series are inserted in (5.8) and (5.9) the features above
are conserved.

From the general constraints (see Appendix)
\begin{eqnarray}
&&\left(A\pm D\right)\left(A'\pm D'\right)=\left(A'\pm
D'\right)\left(A\pm D\right),
\end{eqnarray}
where $A=A\left(\mathbf{x}\right)$,
$A'=A\left(\mathbf{x}'\right)$, etc. it already follows that the
$\mathbf{X}$ (resp. $\mathbf{Y}$) must commute. Thus with
$\mathbf{X}_{(r-2p,2p)}\equiv \mathbf{X}_{(p)}$,
$\mathbf{Y}_{(r-2p-1,2p+1)}\equiv \mathbf{Y}_{(p)}$ one must have,
for all $\left(p,q\right)$
\begin{equation}
\left[\mathbf{X}_{(p)},\mathbf{X}_{(q)}\right]=
\left[\mathbf{Y}_{(p)},\mathbf{Y}_{(q)}\right]=0.
\end{equation}
But our recursion relations imply stronger constraints. From (5.8)
\begin{eqnarray}
&&\mathbf{X}_{(r+1,2p)}=\frac 12
\begin{vmatrix}
  I_{(r)} & K_{(r)} \\
  K_{(r)} & I_{(r)} \\
\end{vmatrix}
\otimes \mathbf{X}_{(r,2p)}+ \frac 12
\begin{vmatrix}
  I_{(r)} & K_{(r)} \\
  -K_{(r)} & -I_{(r)} \\
\end{vmatrix}
\otimes \mathbf{Y}_{(r,2p-1)}
\end{eqnarray}
Note that
\begin{eqnarray}
&&\begin{vmatrix}
  I_{(r)} & K_{(r)} \\
  K_{(r)} & I_{(r)} \\
\end{vmatrix}
\begin{vmatrix}
  I_{(r)} & K_{(r)} \\
  -K_{(r)} & -I_{(r)} \\
\end{vmatrix}=0.
\end{eqnarray}
Along with recursion relations, systematically exploiting the
constraints (A.2-12) of Appendix, one obtains that not only do
$\mathbf{X}_{(p)}$, $\mathbf{X}_{(q)}$ commute for each $r$ but
\begin{equation}
\mathbf{X}_{(p)}\mathbf{X}_{(q)}=\mathbf{X}_{(q)}\mathbf{X}_{(p)}=0,\qquad
p\neq q.
\end{equation}
Analogously one can show
\begin{equation}
\mathbf{Y}_{(p)}\mathbf{Y}_{(q)}=\mathbf{Y}_{(q)}\mathbf{Y}_{(p)}=0,\qquad
p\neq q.
\end{equation}
These are indeed sufficient and necessary conditions for the
eigenvalue spectrum derived for $r=1,2,3,4$ in sec. 4.

Let $\left|\mathbf{v}_p\right\rangle$ denote an eigenstate of
$\mathbf{X}_{(p)}$ with the eigenvalue $\mathbf{X}_{(p)}\left|
\mathbf{v}_p\right\rangle=v_p\left| \mathbf{v}_p\right\rangle$
$\left(v_p\neq 0\right)$. Then for $q\neq p$,
\begin{equation}
\mathbf{X}_{(q)}\left| \mathbf{v}_p\right\rangle=v_p^{-1}
\mathbf{X}_{(q)}\mathbf{X}_{(p)}\left|\mathbf{v}_p\right\rangle=0.
\end{equation}
Hence
\begin{equation}
\left(A+D\right)_r\left| \mathbf{v}_p\right\rangle=
\left(1+\mathbf{x}\right)^{r-2p}
\left(1-\mathbf{x}\right)^{2p}v_p\left|\mathbf{v}_p\right\rangle=0,\qquad
p=0,2,\ldots.
\end{equation}
But there is still one more class of constraints. Starting with
(see (5.2))
\begin{equation}
\hbox{Tr}\left(\mathbf{T}_1\right)=\hbox{Tr}\left(\left(A+D\right)_1\right)=2
\left(1+\mathbf{x}\right)
\end{equation}
and noting that (see (5.8))
\begin{equation}
\hbox{Tr}\left(\mathbf{T}_{r+1}\right)=\hbox{Tr}\left(\left(A+D\right)_{r+1}\right)=
\left(1+\mathbf{x}\right)\hbox{Tr}\left(\left(A+D\right)_{r}\right),
\end{equation}
one obtains
\begin{equation}
\hbox{Tr}\left(\mathbf{T}_{r}\right)=\hbox{Tr}\left(\left(A+D\right)_{r}\right)=
2\left(1+\mathbf{x}\right)^r.
\end{equation}
Hence
\begin{equation}
\sum_{p\neq 0}v_p=0.
\end{equation}
How is this constraint implemented for each $r$? In our examples
(sec. 4, $r=1,2,3,4$) we saw that
\begin{enumerate}
    \item $v_0$ has a multiplicity 2, saturating (5.27).

    \item $v_p$ $\left(p\neq 0\right)$ comes with multiplicity,
    each case providing a subset of zero sum. For each $p$ one
    has, one or more, subsets
\begin{equation}
\sum_{i}v_p^{(i)}=0.
\end{equation}

    \item Here $r$-th roots of unity play a crucial role.
    Typically in (5.29) one has
\begin{equation}
\left(1+\mathbf{x}\right)^{r-2p}\left(1-\mathbf{x}\right)^{2p}\left(\sum_{k=0}^{r-1}e^{\mathbf{i}\frac{2\pi}rk}\right)=0.
\end{equation}

    \item For $r=3$ one has only cube roots of unity. For $r=4$
    one has both 2-plets and 4-plets, (square roots of unity being
    also fourth roots).
\end{enumerate}

Let us now consider possible submultiplets from a more general
point of view. The even and odd subspaces introduced in sec. 4
(even and odd multiplicities of the index 1 distinguishing them)
can be generalized to all $r$, each one closed under the action of
$\left(A+D\right)_r$ and of dimension $2^{r-1}$. In each one there
is exactly 1 state with eigenvalue $\left(1+\mathbf{x}\right)^r$
saturating (5.27). Hence one can now consider separately two base
spaces of dimension $\left(2^{r-1}-1\right)$. When $r$ is a prime
number (say $L$) there is a relative simplicity concerning the
multiplet structure. A Theorem of Fermat (see Ref. 1, Appendix B:
Encounter with a theorem of Fermat) adapted to our case assures
\begin{equation}
2^{L-1}-1=l\cdot L,
\end{equation} where $l$ is an integer. Thus for $L=3,\,5,\,7,\,11$, etc.,
$l=1,3,9,93$, etc.. Hence an integer number of $L$-plets can span
adequately the $2^{L-1}-1$ dimensional space with
$\sum_{k=0}^{r-1}e^{\mathbf{i}\frac{2\pi}{r}k}=0$. When $r$ is not
a prime number each prime factor of $r$ ($4=2\times 2$, $6=2\times
3$, etc.) can lead to submultiplets with zero sum. Finally, if a
singlet occurs in $\left\{e\right\}$ the even subspace (i.e. apart
from $\left(1+\mathbf{x}\right)^r$) it must occur in
$\left\{o\right\}$ the odd one with an opposite sign (ex:
$\left(1-\mathbf{x}\right)^4$ in $\left\{e\right\}$ and
$-\left(1-\mathbf{x}\right)^4$ in $\left\{o\right\}$ for $r=4$).
The number of possibilities increase rapidly with $r$. Our study
remains incomplete concerning the precise number of multiplets and
the multiplicity for each higher $r$. We have however delineated
completely, for all $r$ the dependence of the eigenvalue spectrum
on $\mathbf{x}$ (or the spectral parameter $\theta$).

Leu us note one point. We changed over from $\hat{R}\left(
\theta\right)$ to $\hat{R}\left(\mathbf{x}\right)$ in sec. 3. Thus
gives conveniently a single parameter $\mathbf{x}$ on the
anti-diagonal. With the original normalization
\begin{equation}
\left(1+\mathbf{x}\right)^{r-2p}
\left(1-\mathbf{x}\right)^{2p}\approx
e^{\left(r-2p\right)m_{11}^{(+)}\theta+2pm_{11}^{(-)}\theta}.
\end{equation}

Finally, starting with $\mathbf{T}_1$ or $\mathbf{T}_2$ (of (4.8))
and implementing recursions one can show that the sum of the
elements in each row (and each column) of $T_r$ is
$\left(1+\mathbf{x}\right)^{r}$. The sum of basic components of
$\left\{e\right\}$ and $\left\{o\right\}$
($\left|e_1\right\rangle$, $\left|o_1\right\rangle$ for $r=3,4$ in
sec. 4 and their direct generalizations) thus each corresponds to
the eigenvalue $\left(1+\mathbf{x}\right)^{r}$. The sum of these
two furnishes the total trace of $\mathbf{T}_r$, namely
$2\left(1+\mathbf{x}\right)^r$. All the other states together
contribute zero trace. Moreover, since for our choice of domains
(sec. 3) always $0<\mathbf{x}<1$, $\left(1+\mathbf{x}\right)>1$,
$\left(1-\mathbf{x}\right)<1$ assuming roots of unity in (5.24)
for $v_p$. $\left(1+\mathbf{x}\right)^r$ is the largest
eigenvalue. This is significant in statistical models.

In this context one should note that the special status of
$\mathbf{X}_{(r,0)}$ in the iterative structure. From (5.19) one
obtains the matrices structure
\begin{equation}
2\mathbf{X}_{(r+1,0)}=\begin{vmatrix}
  \mathbf{X}_{(r,0)} & K_{(r)}\mathbf{X}_{(r,0)} \\
  K_{(r)}\mathbf{X}_{(r,0)} & \mathbf{X}_{(r,0)} \\
\end{vmatrix}
\end{equation}
Thus, in each subspace, one can construct iteratively exclusively
$\mathbf{X}_{(r,0)}$ staring from $\mathbf{X}_{(1,0)}$ and the
corresponding $2^{r-1}$ mutually orthogonal eigenstates:
\begin{enumerate}
    \item One with eigenvalue $\left(1+\mathbf{x}\right)^r$.
    \item $2^{r-1}-1$ with eigenvalue zero.
\end{enumerate}
The latter provide non-zero eigenvalues for $\mathbf{X}_{(p)}$
with $p$ non zero. One thus obtains the complete basis of
eigenstates.

\section{Generalizations $\left(n\geq 2\right)$}
\setcounter{equation}{0}

One has for all $n$
\begin{eqnarray}
&&\mathrm{P}P_{ij}^{(\epsilon)}=\frac
12\left\{\left(ji\right)\otimes
\left(ij\right)+\left(\bar{j}\bar{i}\right)\otimes
\left(\bar{i}\bar{j}\right)+\epsilon\left[\left(j\bar{i}\right)\otimes
\left(i\bar{j}\right)+\left(\bar{j}i\right)\otimes
\left(\bar{i}j\right)\right]\right\},\nonumber\\
&&\mathrm{P}P_{i\bar{j}}^{(\epsilon)}=\frac
12\left\{\left(\bar{j}i\right)\otimes
\left(i\bar{j}\right)+\left(j\bar{i}\right)\otimes
\left(\bar{i}j\right)+\epsilon\left[\left(\bar{j}\bar{i}\right)\otimes
\left(ij\right)+\left(ji\right)\otimes
\left(\bar{i}\bar{j}\right)\right]\right\}.
\end{eqnarray}
These lead (for fundamental blocks or $r=1$) to
\begin{eqnarray}
&&\mathbf{T}_{ij}=a_{ji}^{(+)}\left(ji\right)+a_{ji}^{(-)}\left(\bar{j}\bar{i}\right),\qquad
\mathbf{T}_{\bar{i}\bar{j}}=a_{ji}^{(+)}\left(\bar{j}\bar{i}\right)+a_{ji}^{(-)}\left(ji\right),\nonumber\\
&&\mathbf{T}_{i\bar{j}}=a_{ji}^{(-)}\left(j\bar{i}\right)+a_{ji}^{(+)}\left(\bar{j}i\right),\qquad
\mathbf{T}_{\bar{i}j}=a_{ji}^{(-)}\left(\bar{j}i\right)+a_{ji}^{(+)}\left(j\bar{i}\right),
\end{eqnarray}
where $a_{ij}^{(\pm)}=\frac 12\left(e^{m_{ij}^{(+)}\theta}\pm
e^{m_{ij}^{(-)}\theta}\right)$. For $m_{ij}^{(+)}> m_{ij}^{(-)}$
(resp. $m_{ij}^{(+)}< m_{ij}^{(-)}$) all Boltzmann weights are
nonnegative for $\theta>0$ (resp. $\theta<0$). Note that the
fundamental blocks are now $\left(2n\times 2n\right)$ matrices
with only two non-zero elements each. This number does not change
with $n$. In a compact notation with indices $a\in\left\{1,2,
\ldots,2n\right\}$ and $\bar{a}=\left\{2n,2n-1,\ldots,1\right\}$
correspondingly, one can write, for $r=1$,
\begin{equation}
\mathbf{T}_{ab}=a^{(+)}_{ba}\left(ba\right)+a^{(-)}_{ba}\left(\bar{b}\bar{a}\right)
\end{equation}
with (not only $a_{ij}^{(\pm)}=a_{i\bar{j}}^{(\pm)}$, but)
\begin{equation}
a^{(\pm)}_{ba}=a^{(\pm)}_{\bar{b}\bar{a}}=a^{(\pm)}_{b\bar{a}}=a^{(\pm)}_{\bar{b}a}.
\end{equation}
The coproduct rules gives the iterative structure
\begin{equation}
\mathbf{T}_{ab}^{(r+1)}=\sum_c\mathbf{T}_{ac}\otimes
\mathbf{T}_{cb}^{(r)}=
\sum_c\left(a_{ca}^{(+)}\left(ca\right)+a_{ca}^{(-)}\left(\bar{c}\bar{a}\right)\right)\otimes
\mathbf{T}_{cb}^{(r)}.
\end{equation}
The fact that, as in (6.3), only diagonal blocks have diagonal
elements can be easily shown to lead to
\begin{equation}
\Im_r\equiv
\hbox{Tr}\left(\mathbf{T}^{(r)}\right)=2\left(\sum_{i=1}^n
e^{rm_{ii}^{(+)}\theta}\right).
\end{equation}
This is the direct multiparametric generalization of the $4\times
4$ case
\begin{equation}
\Im_r=2e^{rm_{11}^{(+)}\theta}.
\end{equation}

It is instructive to study the case $n=2$ explicitly. Denoting (as
in Ref. 1 with $\bar{1}=2$, $\bar{2}=1$)
\begin{eqnarray}
&&a_{(\pm)}=\frac 12\left(e^{m_{11}^{(+)}\theta}\pm
e^{m_{11}^{(-)}\theta}\right),\qquad d_{(\pm)}=\frac
12\left(e^{m_{22}^{(+)}\theta}\pm
e^{m_{22}^{(-)}\theta}\right),\nonumber\\
&&b_{(\pm)}=\frac 12\left(e^{m_{12}^{(+)}\theta}\pm
e^{m_{12}^{(-)}\theta}\right),\qquad c_{(\pm)}=\frac
12\left(e^{m_{21}^{(+)}\theta}\pm e^{m_{21}^{(-)}\theta}\right)
\end{eqnarray}
for $r=1$ (with now $i=1,2$, $\bar{i}=4,3$ below - a change of
notation convenient for displaying symmetries) one has
\begin{eqnarray}
&&\mathbf{T}_{11}=\begin{vmatrix}
  a_+ & 0 & 0 & 0 \\
  0 & 0 & 0 & 0 \\
  0 & 0 & 0 & 0 \\
  0 & 0 & 0 & a_- \\
\end{vmatrix},\qquad \mathbf{T}_{22}=\begin{vmatrix}
  0 & 0 & 0 & 0 \\
  0 & d_+ & 0 & 0 \\
  0 & 0 & d_- & 0 \\
  0 & 0 & 0 & 0 \\
\end{vmatrix},\nonumber\\
&&\mathbf{T}_{1\bar{1}}=\begin{vmatrix}
  0 & 0 & 0 & a_- \\
  0 & 0 & 0 & 0 \\
  0 & 0 & 0 & 0 \\
  a_+ & 0 & 0 & 0 \\
\end{vmatrix},\qquad \mathbf{T}_{2\bar{2}}=\begin{vmatrix}
  0 & 0 & 0 & 0 \\
  0 & 0 & d_- & 0 \\
  0 & d_+ & 0 & 0 \\
  0 & 0 & 0 & 0 \\
\end{vmatrix},\nonumber\\
&&\mathbf{T}_{12}=\begin{vmatrix}
  0 & 0 & 0 & 0 \\
  c_+ & 0 & 0 & 0 \\
  0 & 0 & 0 & c_- \\
  0 & 0 & 0 & 0 \\
\end{vmatrix},\qquad \mathbf{T}_{1\bar{2}}=\begin{vmatrix}
  0 & 0 & 0 & 0 \\
  0 & 0 & 0 & c_- \\
  c_+ & 0 & 0 & 0 \\
  0 & 0 & 0 & 0 \\
\end{vmatrix},\nonumber
\end{eqnarray}
\begin{eqnarray}
&&T_{21}=\begin{vmatrix}
  0 & b_+ & 0 & 0 \\
  0 & 0 & 0 & 0 \\
  0 & 0 & 0 & 0 \\
  0 & 0 & b_- & 0 \\
\end{vmatrix},\qquad \mathbf{T}_{2\bar{1}}=\begin{vmatrix}
  0 & 0 & b_- & 0 \\
  0 & 0 & 0 & 0 \\
  0 & 0 & 0 & 0 \\
  0 & b_+ & 0 & 0 \\
\end{vmatrix}.
\end{eqnarray}
The remaining 8 blocks are given by
\begin{equation}
\left(a_+,b_+,c_+,d_+\right)\rightleftharpoons
\left(a_-,b_-,c_-,d_-\right)\qquad  \Longrightarrow\qquad
T_{ab}\rightleftharpoons T_{\bar{a}\bar{b}}.\end{equation} One has
\begin{equation}
\Im_1=\hbox{Tr}\left(\mathbf{T}_{11}+\mathbf{T}_{22}+\mathbf{T}_{\bar{1}\bar{1}}+\mathbf{T}_{\bar{2}\bar{2}}\right)
=2\left(\left(a_++a_-\right)+\left(d_++d_-\right)\right)=2
\left(e^{m_{11}^{(+)}\theta}+e^{m_{22}^{(+)}\theta}\right).
\end{equation}
The recursion
\begin{equation}
\Im_{r+1}=\hbox{Tr}\left(\left(a_++a_-\right)\left(\mathbf{T}_{11}^{(r)}+
\mathbf{T}_{\bar{1}\bar{1}}^{(r)}\right)+\left(d_++d_-\right)
\left(\mathbf{T}_{22}^{(r)}+T_{\bar{2}\bar{2}}^{(r)}\right)\right)
\end{equation}
leads to
\begin{equation}
\Im_r=2\left(e^{rm_{11}^{(+)}\theta}+
e^{rm_{22}^{(+)}\theta}\right)
\end{equation}
a particular case of (6.6). For $n=2$ the diagonal blocks, for
example, have the iterative structures
\begin{eqnarray}
&&\mathbf{T}_{11}^{(r+1)}=\begin{vmatrix}
  a_+\mathbf{T}_{11}^{(r)} & 0 & 0 & a_-\mathbf{T}_{\bar{1}1}^{(r)} \\
  c_+\mathbf{T}_{21}^{(r)} & 0 & 0 & c_-\mathbf{T}_{\bar{2}1}^{(r)} \\
  c_+\mathbf{T}_{\bar{2}1}^{(r)} & 0 & 0 & c_-\mathbf{T}_{21}^{(r)} \\
  a_+\mathbf{T}_{\bar{1}1}^{(r)} & 0 & 0 & a_-\mathbf{T}_{11}^{(r)} \\
\end{vmatrix}\qquad \mathbf{T}_{22}^{(r+1)}=\begin{vmatrix}
  0 & b_+\mathbf{T}_{12}^{(r)} & b_-\mathbf{T}_{\bar{1}2}^{(r)} & 0\\
  0 & d_+\mathbf{T}_{22}^{(r)} & d_-\mathbf{T}_{\bar{2}2}^{(r)} & 0 \\
  0 & d_+\mathbf{T}_{\bar{2}2}^{(r)} & d_-\mathbf{T}_{22}^{(r)} & 0 \\
  0 & b_+\mathbf{T}_{\bar{1}2}^{(r)} & b_-\mathbf{T}_{12}^{(r)} & 0 \\
\end{vmatrix}
\end{eqnarray}
$\mathbf{T}_{\bar{1}\bar{1}}^{(r+1)}$ and $\mathbf{T}_{\bar{2}
\bar{2}}^{(r+1)}$ are now obtained by setting respectively in
$\mathbf{T}_{11}^{(r+1)}$ and $\mathbf{T}_{22}^{(r+1)}$
$a_-\mathbf{T}_{\bar{1}\bar{1}}^{(r)}$ for
$a_+\mathbf{T}_{11}^{(r)}$ and so on, systematically in an evident
fashion. Their sum gives the transfer matrix of order
$\left(r+1\right)$ exhibiting the iterative structure. For $r=1$
the transfer matrix is directly diagonal for all $n$ giving
directly the eigenvalues. For $n=2$, for example (consistently
with (6.13))
\begin{equation}
\mathbf{T}^{(1)}=\begin{vmatrix}
  e^{m_{11}^{(+)}\theta} & 0 & 0 & 0 \\
  0 & e^{m_{22}^{(+)}\theta} & 0 & 0 \\
  0 & 0 & e^{m_{22}^{(+)}\theta} & 0 \\
  0 & 0 & 0 & e^{m_{11}^{(+)}\theta} \\
\end{vmatrix},
\end{equation}
with evident generalization for $n>2$. For $n=1$ we have
systematically explored the remarkable structure of the transfer
matrix for all $r$ (see (5.14-18)) and consequences for
eigenstates. A parallel study for $n>1$ is beyond the scope of
this paper. Our results in the section already indicate how the
multiparametric aspects start playing on essential role.

In sec. 5, starting with the $2\times 2$ matrix $K$ (5.3-7) at the
level of $r=1$ and implementing tensor products, powerful
recursion relations were obtained. We started by relating
$\left(A,B,C,D\right)$ among themselves. For $N\geq 2$ one can
similarly relate (for a given pair of indices $(i,j)$) the quartet
$\left(\mathbf{T}_{ij},\mathbf{T}_{\bar{i}\bar{j}},\mathbf{T}_{i\bar{j}},\mathbf{T}_{\bar{i}j}\right)$
given, for $r=1$, by (6.2). Evidently one can relate through
constant matrices only blocks involving the same pair of
parameters $\left(m_{ij}^{(\pm)}\right)$. For this one introduces
the matrix
\begin{equation}\sum_{i=1}^n\left(\left(i\bar{i}\right)+\left(\bar{i}i\right)\right)
\end{equation}
generalizing
$K=\left(1\bar{1}\right)+\left(\bar{1}1\right)=\begin{vmatrix}
  0 & 1 \\
  1 & 0 \\
\end{vmatrix}$. Thus for $n=2$, one has
\begin{equation}
K\otimes K=
\begin{vmatrix}
  0 & 0 & 0 & 1 \\
  0 & 0 & 1 & 0 \\
  0 & 1 & 0 & 0 \\
  1 & 0 & 0 & 0 \\
\end{vmatrix}=K^{(2)}.\end{equation}
The generalization is evident. From (6.8)
\begin{eqnarray}
&&K^{(2)}\left(\mathbf{T}_{11},\mathbf{T}_{\bar{1}\bar{1}},\mathbf{T}_{1\bar{1}},\mathbf{T}_{\bar{1}1}\right)=
\left(\mathbf{T}_{1\bar{1}},\mathbf{T}_{\bar{1}1},\mathbf{T}_{11},\mathbf{T}_{\bar{1}\bar{1}}\right),\nonumber\\
&&K^{(2)}\left(\mathbf{T}_{11},\mathbf{T}_{\bar{1}\bar{1}},\mathbf{T}_{1\bar{1}},\mathbf{T}_{\bar{1}1}\right)K^{(2)}=
\left(\mathbf{T}_{\bar{1}\bar{1}},\mathbf{T}_{11},\mathbf{T}_{\bar{1}1},\mathbf{T}_{1\bar{1}}\right)
\end{eqnarray}
with exactly analogous results for other subsets. For each $n$ one
has, in evident notations, for fixed $\left(i,j\right)$
\begin{eqnarray}
&&K^{(n)}\left(\mathbf{T}_{ij},\mathbf{T}_{\bar{i}\bar{j}},\mathbf{T}_{i\bar{j}},\mathbf{T}_{\bar{i}j}\right)=
\left(\mathbf{T}_{i\bar{j}},\mathbf{T}_{\bar{i}j},\mathbf{T}_{ij},\mathbf{T}_{\bar{i}\bar{j}}\right),\nonumber\\
&&K^{(n)}\left(\mathbf{T}_{ij},\mathbf{T}_{\bar{i}\bar{j}},\mathbf{T}_{i\bar{j}},\mathbf{T}_{\bar{i}j}\right)K^{(n)}=
\left(\mathbf{T}_{\bar{i}\bar{j}},\mathbf{T}_{ij},\mathbf{T}_{\bar{i}j},\mathbf{T}_{i\bar{j}}\right).
\end{eqnarray}
We will not attempt to explore in the present paper the
applications of such relations generalizing our results of sec. 5.

\section{Spin chain Hamiltonians}
\setcounter{equation}{0}

Our construction of odd dimensional Hamiltonians (sec. 4, Ref. 1)
can be adapted to the present even dimensional cases as follows.
One has, taking derivatives and setting $\theta=0$, for $n=1$,
\begin{equation}
\dot{\hat{R}}\left(0\right)=\begin{vmatrix}
  \mathbf{x}_+ & 0 & 0 & \mathbf{x}_- \\
  0 & \mathbf{x}_+ & \mathbf{x}_- & 0 \\
  0 & \mathbf{x}_- & \mathbf{x}_+ & 0 \\
  \mathbf{x}_- & 0 & 0 & \mathbf{x}_+ \\
\end{vmatrix},
\end{equation}
where $\mathbf{x}_{\pm}=\frac 12\left(m_{11}^{(+)}\pm
m_{11}^{(-)}\right)$ and $\dot{\hat{R}}\left(0\right)$ is obtained
by setting $\theta=0$ in $\frac
{d}{d\theta}\hat{R}\left(\theta\right)$. For $n=2$, setting
(starting from (6.8))
\begin{eqnarray}
&&\hat{a}_{\pm}=\frac 12\left(m_{11}^{(+)}\pm
m_{11}^{(-)}\right),\qquad \hat{d}_{\pm}=\frac
12\left(m_{22}^{(+)}\pm m_{22}^{(-)}\right),\nonumber\\
&&\hat{b}_{\pm}=\frac 12\left(m_{12}^{(+)}\pm
m_{12}^{(-)}\right),\qquad \hat{c}_{\pm}=\frac
12\left(m_{21}^{(+)}\pm m_{21}^{(-)}\right)
\end{eqnarray}
from (2.7) and
\begin{equation}
\dot{\hat{R}}\left(0\right)=\begin{vmatrix}
  \hat{D}_{11} & 0 & 0 & \hat{A}_{1\bar{1}} \\
  0 & \hat{D}_{22} & \hat{A}_{2\bar{2}} & 0 \\
  0 & \hat{A}_{\bar{2}2} & \hat{D}_{\bar{2}\bar{2}} & 0 \\
  \hat{A}_{\bar{1}1} & 0 & 0 & \hat{D}_{\bar{1}\bar{1}} \\
\end{vmatrix},
\end{equation}
where
\begin{eqnarray}
&&\hat{D}_{11}=\hat{D}_{\bar{1}\bar{1}}=\left(\hat{a}_+,\hat{b}_+,\hat{b}_+,\hat{a}_+\right)_{\hbox{diag.}}\qquad
\hat{D}_{22}=\hat{D}_{\bar{2}\bar{2}}=\left(\hat{c}_+,\hat{d}_+,\hat{d}_+,\hat{c}_+\right)_{\hbox{diag.}}\nonumber\\
&&\hat{A}_{1\bar{1}}=\hat{A}_{\bar{1}1}=\left(\hat{a}_-,\hat{b}_-,\hat{b}_-,\hat{a}_-\right)_{\hbox{anti-diag.}}\qquad
\hat{A}_{2\bar{2}}=\hat{A}_{\bar{2}2}=\left(\hat{c}_-,\hat{d}_-,\hat{d}_-,\hat{c}_-\right)_{\hbox{anti-diag.}}.
\end{eqnarray}
The extension of our formalism for $n>2$ is straightforward.

For $r$ sites the standard result for the Hamiltonian is (see
sources cited in Ref. 1)
\begin{equation}
\mathrm{H}=\sum_{k=1}^rI\otimes\cdots\otimes
\dot{\hat{R}}\left(0\right)_{k,k+1}\otimes\cdots \otimes I,
\end{equation}
where for circular boundary conditions $k+1=r+1\approx 1$. We
intend to present a more complete study of our spin chain
elsewhere. But here one already sees how the two aspects,
multistate (higher spins at each site) and multiparameter
$\left(m_{ij}^{(\pm)}\right)$ get directly associated for our
hierarchy. In our previous papers \cite{R3,R6} and here again (see
(3.4), (3.5)) we showed how the passage to imaginary parameters
can lead to unitary $\hat{R}\left(\theta\right)$. The
corresponding $\dot{\hat{R}}\left(0\right)$ has only an overall
factor $\mathbf{i}$ which can be extracted from the sum (7.5).

\section{Remarks}
\setcounter{equation}{0}

\paragraph{\bf I. Status of eigenstates for $n=1$:} For the
simplest $4\times 4$ braid matrix in our hierarchy the
construction of the eigenvalues and eigenfuntions of transfer
matrices of successive orders $\left(r=1,\,2,\,3,\,4,\,
\hbox{etc.}\right)$ has attained the following stage:
\begin{description}

\item[(1)] The transfer matrix at the level $r$ has been expressed
in the form of a $2^r\times 2^r$ matrix
 \begin{equation}
 \mathbf{T}_r=\sum_{p=0,1,2,\ldots,p_m}\left(1+\mathbf{x}\right)^{r-2p}\left(1-\mathbf{x}\right)^{2p}\mathbf{X}_{(p)},
\end{equation}
where $\mathbf{x}=\tanh\frac 12\left(m_{11}^{(+)}-
m_{11}^{(-)}\right)\theta$ and $2p_m=r$ (resp. $r-1$) for $r$ even
(resp. odd) ($2p_m=r-\left(1-(-1)^r\right)/2$). $\theta$ being the
spectral parameter, $m_{11}^{(\pm)}$ two free parameters (see
(3.3)). The matrices $\mathbf{X}_{(p)}$ are constant ones
($\mathbf{x}$-independent) and satisfy
\begin{equation}
\mathbf{X}_{(p)}\mathbf{X}_{(q)}=\mathbf{X}_{(q)}\mathbf{X}_{(p)}=0,\qquad
p\neq q.
\end{equation}
They can be computed systematically via the recursion relations
(for $r\longrightarrow r+1$).

\item[(2)] For
\begin{eqnarray}
&&\mathbf{X}_{(p)}\left|p_{(i)}\right\rangle=v_{(p,i)}\left|p_{(i)}\right\rangle,\qquad
\left(v_{(p,i)}\neq 0\right),\nonumber\\
&& \mathbf{X}_{(q)}\left|p_{(i)}\right\rangle=0,\qquad p\neq q
\end{eqnarray}
and $v_{(p,i)}$ denote phase factors which come in multiplets of
zero sum formed by roots of unity corresponding to $r$ and its
prime factors (see examples in sec. 4). The index "$i$" denotes
such possible multiplicity of each $p$. The exception of the zero
sum rule corresponds to $p=0$. One obtains for each $r$, twice
$v_{(0)}=1$ giving
\begin{equation}
\hbox{Tr}\left(\mathbf{T}_r\right)=2\left(1+\mathbf{x}\right)^r
\end{equation}
a general constraint obtained via recursions. This multiplicity 2
corresponds to two $2^{r-1}$ dimensional subspaces (see "even",
"odd" subspaces defined in sec. 4) each providing just one
eigenstate of $\mathbf{X}_{(0)}$ with non-zero eigenvalue
$\left(1+\mathbf{x}\right)^r$.

\item[(3)] Thus the problem has been reduced to construction of
eigenstates of each $\mathbf{X}_{(p)}$ separately, reducing the
dimension by considering each (even, odd) subspace by turn. This
involves solving sets of linear constraints with only positive and
negative integers as coefficients. One finally keeps only the
non-zero eigenvalues for each $p$ , they being associated with
zero eigenvalues for $\mathbf{X}_{(q)}$, $q\neq p$. In fact as
noted below (5.33), it suffices to construct the full sect of
mutually orthogonal eigenstates of $\mathbf{X}_{(r,0)}$, all but
one in each subspace having eigenvalue zero.

\item[(4)] As already stated (sec. 5) our results remain
incomplete concerning the pattern of possible multiplets and
submultiplets corresponding to roots of unity provided by $r$ and
its prime factors - and multiplicities of such multiplets. A
canonical enumeration when $r$ has a very large number of prime
factors seems to be an unlikely possibility. Nor have we
established rigorously that $v_{(p,i)}$ in (8.3) are always $\pm
1$ or higher roots of unity phases factors. This what happens in
examples $\left(r\leq 4\right)$ of sec. 4 and directly leads to
the following obligatory constraint (8.4).

\item[(4)] We have completely, and for all $r$, extracted the
$\theta$-dependence of eigenvalues in (8.1).

\end{description}

\paragraph{\bf II. Comparisons with standard six vertex and eight vertex models:}

From (2.5) and (3.1) our Yang-Baxter matrix (for $n=1$) with
corresponding normalizations, is
\begin{equation}
R\left(\theta\right)=\mathrm{P}\hat{R}\left(\theta\right)=\begin{vmatrix}
  a_+ & 0 & 0 & a_- \\
  0 & a_- & a_+ & 0 \\
  0 & a_+ & a_- & 0 \\
  a_- & 0 & 0 & a_+ \\
\end{vmatrix},
\end{equation}
where $a_{\pm}=\frac 12\left(e^{m_{11}^{(+)}\theta}\pm
e^{m_{11}^{(-)}\theta}\right)$ and, equivalently,
\begin{equation}
R\left(\mathbf{x}\right)=\begin{vmatrix}
  1 & 0 & 0 & \mathbf{x} \\
  0 & \mathbf{x} & 1 & 0 \\
  0 & 1 & \mathbf{x} & 0 \\
  \mathbf{x} & 0 & 0 & 1 \\
\end{vmatrix},
\end{equation}
where $\mathbf{x}=\tanh\frac
12\left(m_{11}^{(+)}-m_{11}^{(-)}\right)\theta$. Let us now
compare this to the very well known $4\times 4$ six vertex and
eight vertex models - concerning which it is sufficient to a cite
a standard text book \cite{R7} and review articles \cite{R5,R8}
which cite basic sources. All such cases are of the form
\begin{equation}
R\left(\theta\right)=\begin{vmatrix}
  a & 0 & 0 & d \\
  0 & b & c & 0 \\
  0 & c & b & 0 \\
  d & 0 & 0 & a \\
\end{vmatrix},
\end{equation}
For (8.6) $a=c=1$, $b=d=\mathbf{x}$ and for (8.5) $a=c=a_+$,
$b=d=a_-$. For six vertex, crucially, $d=0$ and
$\left(a,b,c\right)$ is being given, according to the regime, by
circular or hyperbolic functions. For eight vertex, famously,
elliptic functions appear and $d\neq 0$. For our case $d\neq 0$
but arguably, one has maximal simplicity and symmetry compatible
with non-trivial solution for a $4\times 4$ Yang-Baxter (or braid)
matrix. We have all eight vertices but (with, say (8.5)) ( see
fig. 1)\begin{figure}[ht]
\centerline{\includegraphics[height=2cm]{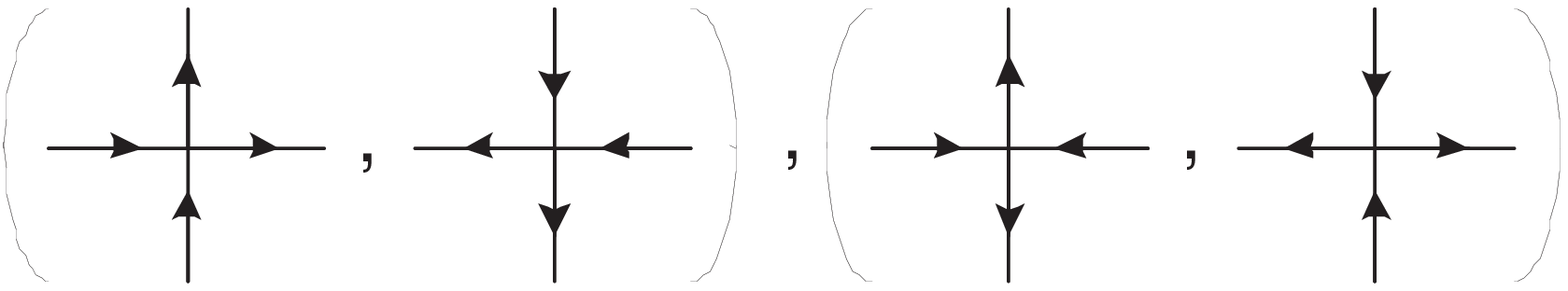}}
\centerline{fig. 1}
\end{figure}
 corresponding to $a=c=a_+$
and (see fig. 2)\begin{figure}[ht]
\centerline{\includegraphics[height=2cm]{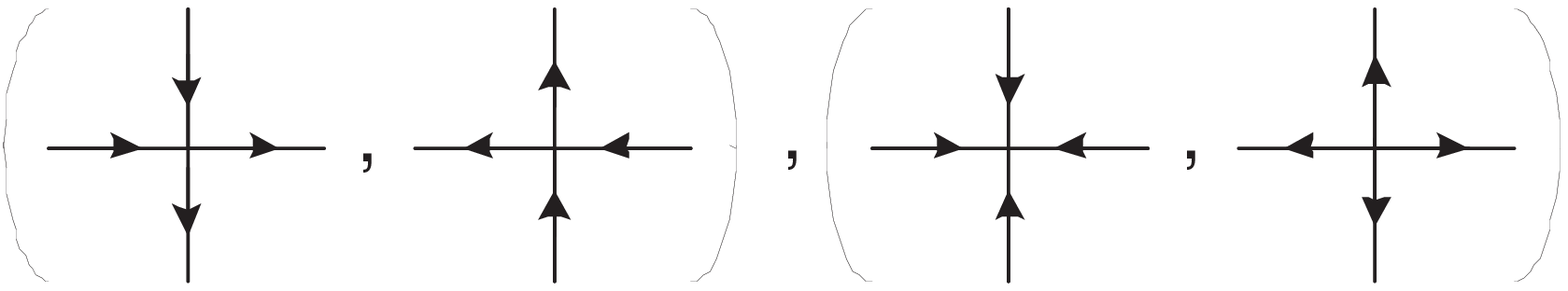}}
\centerline{fig. 2}
\end{figure}
 corresponding to $b=d=a_-$. In six vertex the last two
vertices are excluded $\left(d=0\right)$. We will not study, in
this paper, the implications of the results below (8.7) concerning
various properties of our model (compare the relevant detailed
study of eight vertex in Ref. 7). But we would like to contrast
our approach to the construction of eigenstates and extraction of
eigenvalues of $\mathbf{T}_r$ with that via Bethe ansatz in
standard six vertex models \cite{R5}. The systematic study of
$R\mathbf{TT}$ constraints are particularly relevant (see Appendix
A).

In six vertex the Bethe ansatz construction involves pushing
$\left(A\left(\theta\right)+D\left(\theta\right)\right)$ through
the product $B\left(\theta_1\right)B\left(\theta_2\right)\ldots
B\left(\theta_r\right)$ acting on one single state
$\left|\begin{matrix}1\\0\\\end{matrix}\right\rangle_1\otimes
\left|\begin{matrix}
1\\0\\\end{matrix}\right\rangle_2\otimes\cdots\left|\begin{matrix}
1\\0\\\end{matrix}\right\rangle_r\equiv \left|11\ldots
1\right\rangle_r$ and eliminating unwanted terms to obtain a
complete set of eigenstates corresponding to the sets of resulting
constraints. This involves solving nonlinear equations. Our
resorts to programs and numerical studies for higher $r$'s.

In our case recursion relations for $\left(A\pm D\right)_r$ (for
$r\longrightarrow r+1$) are sufficient to attain the stage
systematically presented in part (I) of this section. One solves ,
at each stage, linear equations with integer coefficients. In fact
since each $r$ (see (5.5-7)) $\left(B\pm C\right)_r=K_{(r)}
\left(A\pm D\right)_r$, where $K_{(r)}$ in the $2^r\times 2^r$
matrix with $2^{r-1}$-times $K=\begin{vmatrix}
  1 & 0 \\
  0 & 1 \\
\end{vmatrix}$ on diagonal , the actions of $\left(B\pm
C\right)_r$ follows immediately from those of $\left(A\pm
D\right)_r$ . They need hardly be studied separately. However that
$K_{(r)}$ and hence $\left(B\pm C\right)_r$ connect the even and
odd subspaces. On the other hand $\left(A+D\right)$ cannot be
pushed through a product of $B\left(\theta\right)$'s displayed
above. The nearest approaches are typically our (A.10) and (A.12)
(changes of signs in $\left(A\pm D\right)$ and $\left(B\pm
C\right)$ are to be noted as $\left(A+D\right)$ is pushed
through). Subtle analytic properties, unlike six vertex and
particularly eight vertex models, play no role in our case. One
has only to look at the $x$ (or $\theta$) dependence in (8.1).

In a simple situation the contrast between our model and the
standard eight vertex one shows up very clearly: For the
asymptotic case $\theta$ tending to infinity one sets
$\mathbf{x}=1$ in (8.6) to obtain
\begin{equation}
R=\hat{R}=\begin{vmatrix}
  1 & 0 & 0 & 1 \\
  0 & 1 & 1 & 0 \\
  0 & 1 & 1 & 0 \\
  1 & 0 & 0 & 1 \\
\end{vmatrix}.
\end{equation}
Whereas, choosing adequately the normalizing factor (see sec. 7 of
\cite{R9} and in particular eq. (7.13)) for the standard case the
corresponding limiting form is
\begin{equation}
R= \left(1,q^{-1},q^{-1},1\right)_{\hbox{diag.}}\neq\hat{R}.
\end{equation}
When diagonalized (or block diagonalized) our classes of matrices,
in general, lose braid (or Yang-Baxter) property unless the
diagonalizer has a corrected tensor structure (sec. 4 and Addendum
of ref. 3).

\paragraph{\bf III. passage to higher dimensions $n\geq 2$:} As
already emphasized before, a major interest of our simple model
for $n=1$ (the $4\times 4$ braid matrix) is that it is the first
one in a hierarchy of $(2n)^2\times (2n)^2$ braid matrices with
$2n^2$ free parameters at each level. Some of the simplicity of
the $n=1$ case is inevitably lost as $n$ increases. But we have
pointed out in sections 6 and 7 how certain basic features vary in
a simple, canonical fashion. Thus for example,
\begin{description}
    \item[(1)]
\begin{equation}
Tr\left(\mathbf{T}^{(r)}\right)=2\left(\sum_{i=1}^n
e^{rm_{ii}^{(+)}\theta}\right),\qquad \hbox{for all $n$}.
\end{equation}
    \item[(2)] Each blocks $\mathbf{T}_{ij}$ of the $\mathbf{T}$ matrix has just two
    non-zero elements (see (6.2) and (6.9)) out of $\left(2n\times
    2n\right)$, namely $a_{ji}^{(\pm)}\left(\theta\right)$ for all
    $n$.
    \item[(3)] For a fixed pair of indices $\left(i,j\right)$ the
    blocks $\left(\mathbf{T}_{ij},\mathbf{T}_{\bar{i}\bar{j}},\mathbf{T}_{i\bar{j}},\mathbf{T}_{\bar{i}j}\right)$ can be quite simply related among
    themselves (see (6.18)) via direct generalizations of the
    matrix $K=\begin{vmatrix}
      0 & 1 \\
      1 & 0 \\
    \end{vmatrix}$ for the $4\times 4$ case. For $n=1$ such
    relations led to recursion relations yielding (8.1-3).

 \item[(4)] Spin chain Hamiltonians present (see eqs.
(7.1-5)) a simple canonical sequence as $n$ increases.
\end{description}

We hope to study the higher dimensional cases more fully
elsewhere.

\vskip 0.5cm

\noindent{\bf Acknowledgments:} {\em One of us (BA) wants to thank
Pierre Collet and Paul Sorba for precious help. }

\begin{appendix}

\section{$R\mathbf{TT}$ constraints $\left(n=1\right)$}
\setcounter{equation}{0}

From (3.10) one obtains with $K$ from (3.9)
\begin{eqnarray}
&&\begin{vmatrix}
  A & B \\
  C & D \\
\end{vmatrix}\otimes\begin{vmatrix}
  A' & B' \\
  C' & D' \\
\end{vmatrix}-\begin{vmatrix}
  A' & B' \\
  C' & D' \\
\end{vmatrix}\otimes\begin{vmatrix}
  A & B \\
  C & D \\
\end{vmatrix}\nonumber\\
&&=\mathbf{x}''\left(\begin{vmatrix}
  B' & A' \\
  D' & C' \\
\end{vmatrix}\otimes\begin{vmatrix}
  B & A \\
  D & C \\
\end{vmatrix}-\begin{vmatrix}
  C & D \\
  A & B \\
\end{vmatrix}\otimes\begin{vmatrix}
  C' & D' \\
  A' & B' \\
\end{vmatrix}\right)
\nonumber\\
&&=\mathbf{x}''\left(\begin{vmatrix}
  B & A \\
  D & C \\
\end{vmatrix}\otimes\begin{vmatrix}
  B' & A' \\
  D' & C' \\
\end{vmatrix}-\begin{vmatrix}
  C' & D' \\
  A' & B' \\
\end{vmatrix}\otimes\begin{vmatrix}
  C & D \\
  A & B \\
\end{vmatrix}\right).
\end{eqnarray}
The last step follows from $\mathbf{x}\rightleftharpoons
\mathbf{x}'$ on both sides of the first two expressions since
under this interchange (see (3.7)) $\mathbf{x}''\longrightarrow
-\mathbf{x}''$. The consistency of the last two implies that each
element of the total matrix on the right (apart from the factor
$\mathbf{x}''$) must be symmetric in
$\left(\mathbf{x},\mathbf{x}'\right)$. This can indeed be verified
starting from $r=1,2,3,\ldots$ using the standard construction in
sec. 3.

From the last two steps (the factor $\mathbf{x}''$ cancelling) one
obtains, with only upper or lower signs,
\begin{eqnarray}
&&\left(A\pm D\right)\left(A'\pm D'\right)=\left(A'\pm
D'\right)\left(A\pm D\right),\nonumber\\
&&\left(B\pm C\right)\left(B'\pm C'\right)=\left(B'\pm
C'\right)\left(B\pm C\right),\nonumber\\
&&\left(A\pm D\right)\left(B'\pm C'\right)=\left(A'\pm
D'\right)\left(B\pm C\right),\nonumber\\
&& \left(B\pm C\right)\left(A'\pm D'\right)=\left(B'\pm
C'\right)\left(A\pm D\right).
\end{eqnarray}
One gets 8 relations of the type
\begin{eqnarray}
&&\left(M_1M_2'-M_1'M_2\right)=\mathbf{x}''\left(M_3'M_4-M_5M_6'\right)
=\mathbf{x}''\left(M_3M_4'-M_5'M_6\right),
\end{eqnarray}
with
\begin{eqnarray}
&&\left(M_1,M_2,M_3,M_4,M_5,M_6\right)\nonumber\\
&&=\left\{
\left(A,B,B,A,C,D\right),\left(A,C,B,D,C,A\right),\left(B,A,A,B,D,C\right)
\left(B,D,A,C,D,B\right),\right.\nonumber\\
&&\left.\left(C,A,D,B,A,C\right),
\left(C,D,D,C,A,B\right),\left(D,B,C,A,B,D\right),\left(D,C,C,D,B,A\right)
\right\},
\end{eqnarray}
In sec. 5 (see (5.4-7)) we obtained, for any $\mathbf{x}$ and all
$r$,
\begin{equation}
\left(B,C\right)=K_{(r)}\left(A,D\right),\qquad
\left(A,D\right)=K_{(r)}\left(B,C\right),
\end{equation}
where ($\mathbf{x}$-independent ) matrix $K_{(r)}$ is given. Hence
multiplying in (A.3) by $K_{(r)}$ on the left one gets another set
with $\left(M_1,M_3,M_5\right)$ replaced by
$K_{(r)}\left(M_1,M_3,M_5\right)$, where thus
\begin{equation}
\left(A,B,B,A,C,D\right)\longrightarrow\left(B,B,A,A,D,D\right)
\end{equation}
and so on. For another class of relations we introduce
\begin{equation}
\mathbf{f}^{(\pm)}=\frac 12\left(\mathbf{x}''\pm \frac
1{\mathbf{x}''}\right)
\end{equation}
From (3.7)
\begin{equation}
{\bf f}^{(+)}=\coth\mu\left(\theta-\theta'\right),\qquad {\bf
f}^{(-)}=-\hbox{cosech}\mu\left(\theta-\theta'\right)
\end{equation}
One can show that (with upper or lower signs)
\begin{eqnarray}
&&\left(A\pm D\right)\left(A'\mp D'\right)={\bf f}^{(+)}\left(B\pm
C\right)\left(B'\mp C'\right)+{\bf f}^{(-)}\left(B'\pm
C'\right)\left(B\mp C\right),\\
&&\left(A\pm D\right)\left(B'\mp C'\right)={\bf f}^{(+)}\left(B\pm
C\right)\left(A'\mp D'\right)+{\bf f}^{(-)}\left(B'\pm
C'\right)\left(A\mp D\right).
\end{eqnarray}
Again, as for (A.5), (A.6), multiplying the above from the left by
$K_{(r)}$ one gets another set of relations such that
\begin{eqnarray}
&&\left(B\pm C\right)\left(B'\pm C'\right)={\bf f}^{(+)}\left(A\pm
D\right)\left(A'\mp D'\right)+{\bf f}^{(-)}\left(A'\pm
D'\right)\left(A\mp D\right).
\end{eqnarray}
and so on. For "two-steep" relations redefine $x''$ with indices
$\mathbf{x}_{ij}=\frac
{\mathbf{x}_i-\mathbf{x}_j}{1-\mathbf{x}_i\mathbf{x}_j}$ and ${\bf
f}^{(+)}_{ij}=\coth\mu\left(\theta_i-\theta_j\right)$, ${\bf
f}^{(-)}_{ij}=$-cosech$\mu\left(\theta_i-\theta_j\right)$
correspondingly. One obtains, for example, for arguments
$\left(\mathbf{x}_1,\mathbf{x}_2,\mathbf{x}_3\right)$
\begin{eqnarray}
&&\left(A+D\right)_1\left(B-C\right)_2\left(B+C\right)_3={\bf
f}^{(+)}_{12}{\bf f}^{(+)}_{23}\left(B+C\right)_1\left(B-
C\right)_2\left(A+D\right)_3+\nonumber\\
&&\phantom{\left(A+D\right)_1\left(B-C\right)_2\left(B+C\right)_3=}{\bf
f}^{(+)}_{12}{\bf f}^{(-)}_{23}\left(B+C\right)_1\left(B-
C\right)_3\left(A+D\right)_2+\nonumber\\
&&\phantom{\left(A+D\right)_1\left(B-C\right)_2\left(B+C\right)_3=}{\bf
f}^{(-)}_{12}{\bf f}^{(+)}_{13}\left(B+C\right)_2\left(B-
C\right)_1\left(A+D\right)_3+\nonumber\\
&&\phantom{\left(A+D\right)_1\left(B-C\right)_2\left(B+C\right)_3=}{\bf
f}^{(-)}_{12}{\bf f}^{(-)}_{13}\left(B+C\right)_2\left(B-
C\right)_3\left(A+D\right)_1.
\end{eqnarray}
See the relevant remarks in sec. 8.

\end{appendix}


\vskip 0.5cm

\begin{thebibliography}{99}

\bibitem{R1} B. Abdesselam and A. Chakrabarti, {\it A nested sequence of projectors: (2) Multiparameter
multistate statistical models, Hamiltonians, $S$-matrices}, Jour.
Math. Phys. {\bf 47}(2006) 053508.

\bibitem{R2} A. Chakrabarti, {\it A nested sequence of projectors and
corresponding braid matrices $\hat R(\theta)$: (1) Odd
dimensions}, Jour. Math. Phys. {\bf 46} (2005) 063508.

\bibitem{R3} B. Abdesselam, A. Chakrabarti, V.K. Dobrev and S.G. Mihov,
{\it Higher dimensional multiparameter unitary and nonunitary
braid matrices: even dimensions}, Jour. Math. Phys. {\bf 48}
(2007) 103505.

\bibitem{R4} D. Arnaudon, A. Chakrabarti, V.K. Dobrev and S.G. Mihov,
{\it Exotic Bialgebra S\O3: Representations, Baxterisation and
Applications}, Ann. H. poincar\'e {\bf 7} (2006) 1351.

\bibitem{R5} H.J. De Vega, {\it Yang-Baxter algebras, integrable theories and quantum groups},
Int. Jour. Mod. Phys. A vol. {\bf 4} (1989) 2371.

\bibitem{R6} B. Abdesselam, A. Chakrabarti, V.K. Dobrev and S.G. Mihov,
{\it Higher Dimensional Unitary Braid Matrices: Construction,
Associated Structures and Entanglements}, Jour. Math. Phys. {\bf
48} (2007) 053508.

\bibitem{R7} R.J. Baxter, {\it Exactly solved models in statistical
mechanics}, Acad. Press (1982).

\bibitem{R8} H. Saleur and J.B. Zuber, {\it Integrable lattice models
and quantum groups}, in the proceedings of the 1990 Trieste Spring
School on String Theory and Quantum Gravity.

\bibitem{R9} A. Chakrabarti, {\it Canonical factorization and diagonalization of Baxterized braid matrices:
Explicit constructions and applications}, Jour. Math. Phys. {\bf
44} (2003) 5320.


\end{thebibliography}
\end{document}